# Near optimal thresholding estimation of a Poisson intensity on the real line

Patricia Reynaud-Bouret[1] and Vincent Rivoirard[2]

**Abstract** The purpose of this paper is to estimate the intensity of a Poisson process $N$ by using thresholding rules. In this paper, the intensity, defined as the derivative of the mean measure of $N$ with respect to $ndx$ where $n$ is a fixed parameter, is assumed to be non-compactly supported. The estimator $\tilde{f}_{n,\gamma}$ based on random thresholds is proved to achieve the same performance as the oracle estimator up to a possible logarithmic term. Then, minimax properties of $\tilde{f}_{n,\gamma}$ on Besov spaces $\mathcal{B}^{\alpha}_{p,q}$ are established. Under mild assumptions, we prove that

$$\sup_{f \in \mathcal{B}^{\alpha}_{p,q} \cap \mathbb{L}_\infty} \mathbb{E}(\|\tilde{f}_{n,\gamma} - f\|_2^2) \leq C \left(\frac{\log n}{n}\right)^{\frac{\alpha}{\alpha + \frac{1}{2} + \left(\frac{1}{2} - \frac{1}{p}\right)_+}}$$

and the lower bound of the minimax risk for $\mathcal{B}^{\alpha}_{p,q} \cap \mathbb{L}_\infty$ coincides with the previous upper bound up to the logarithmic term. This new result has two consequences. First, it establishes that the minimax rate of Besov spaces $\mathcal{B}^{\alpha}_{p,q}$ with $p \leq 2$ when non compactly supported functions are considered is the same as for compactly supported functions up to a logarithmic term. When $p > 2$, the rate exponent, which depends on $p$, deteriorates when $p$ increases, which means that the support plays a harmful role in this case. Furthermore, $\tilde{f}_{n,\gamma}$ is adaptive minimax up to a logarithmic term.

## 1 Introduction

The goal of the present paper is to derive a data-driven thresholding method to estimate the intensity of a Poisson process on the real line.

Poisson processes have been used for years to model a wide variety of situations, and in particular data whose maximal size is a priori unknown. For instance, in finance, Merton [29] introduces Poisson processes to model stock-price changes of extraordinary magnitude. In geology, Uhler and Bradley [32] use Poisson processes to model the occurrences of petroleum reservoirs whose size is highly inhomogeneous. Actually, if we only focus on the size of the jumps in Merton's model or on the sizes of individual oil reservoirs, these models consist in an inhomogeneous Poisson process with heavy-tailed intensities (see [19] for a precise formalism for the financial example). So, our goal is to provide data-driven estimation of a Poisson intensity with as few support assumptions as possible.

---

[1]CNRS and Département de Mathématiques et Applications, ENS-Paris, 45 Rue d'Ulm, 75230 Paris Cedex 05, France. Email: reynaud@dma.ens.fr

[2]Equipe Probabilité, Modélisation et Statistique, Laboratoire de Mathématique, CNRS UMR 8628, Université Paris Sud, 91405 Orsay Cedex, France. Département de Mathématiques et Applications, ENS-Paris, 45 Rue d'Ulm, 75230 Paris Cedex 05, France. Email: Vincent.Rivoirard@math.u-psud.fr

Of course, many adaptive methods have been proposed to deal with Poisson intensity estimation. For instance, Rudemo [31] studied data-driven histogram and kernel estimates based on the cross-validation method. Donoho [15] fitted the universal thresholding procedure proposed by Donoho and Johnstone [17] by using the Anscombe's transform. Kolaczyk [28] refined this idea by investigating the tails of the distribution of the noisy wavelet coefficients of the intensity. For a particular inverse problem, Cavalier and Koo [10] first derived optimal estimates in the minimax setting. More precisely, for their tomographic problem, Cavalier and Koo [10] pointed out minimax thresholding rules on Besov balls. By using model selection, other optimal estimators have been proposed by Reynaud-Bouret [30] or Willet and Nowak [33].

To derive sharp theoretical results, these methods need to assume that the intensity has a known bounded support and belongs to $\mathbb{L}_\infty$. Model selection may allow to remove the assumption on the support. See oracle results established by [19] who nevertheless assumes that the intensity belongs to $\mathbb{L}_\infty$. We have to mention that the model selection methodology proposed by Baraud and Birgé [7], [4] is "assumption-free" as well. However, as explained by Birgé [7], it is too computationally intensive to be implemented. Besides, in [7], [4] and [19], minimax performances on classical functional spaces are derived only for compactly supported signals.

In the present paper, to estimate the intensity of a Poisson process, we propose an easily implementable thresholding rule specified in the next section. This procedure is near optimal under oracle and minimax points of view. We do not assume that the support of the intensity is known or even finite and most of the time, the signal to estimate may be unbounded.

## 1.1   The thresholding procedure and main result

In the sequel, we consider a Poisson process on the real line, denoted $N$, whose mean measure $\mu$ is finite and absolutely continuous with respect to the Lebesgue measure (see Section 2.1 where we recall classical facts on Poisson processes). Given $n$ a positive integer, we introduce $f \in \mathbb{L}_1(\mathbb{R})$ the intensity of $N$ as

$$f(x) = \frac{d\mu_x}{ndx}.$$

Since $f$ belongs to $\mathbb{L}_1(\mathbb{R})$, the total number of points of the process $N$, denoted $N_\mathbb{R}$, satisfies $\mathbb{E}(N_\mathbb{R}) = n\|f\|_1$ and $N_\mathbb{R} < \infty$ almost surely. In the sequel, $f$ will be held fixed and $n$ will go to $+\infty$. The introduction of $n$ could seem artificial, but it allows to present our asymptotic theoretical results in a meaningful way. In addition, our framework is equivalent to the observation of a $n$-sample of a Poisson process with common intensity $f$ with respect to the Lebesgue measure. Since $N$ is a random countable set of points, we denote by $dN$ the discrete random measure $\sum_{T \in N} \delta_T$. Hence we have for any compactly supported function $g$, $\int g(x) dN_x = \sum_{T \in N} g(T)$. Now, our goal is to estimate $f$ by using the realizations of $N$.

For this purpose, we assume that $f$ belongs to $\mathbb{L}_2(\mathbb{R})$ and we use the decomposition of $f$ on one of the biorthogonal wavelet bases described in Section 2.2. We recall that, as classical orthonormal wavelet bases, biorthogonal wavelet bases are generated by dilations and translations of father and mother wavelets. But considering biorthogonal wavelets allows to distinguish, if necessary, wavelets for analysis (that are piecewise constant functions in this paper) and wavelets for reconstruction with a prescribed number of continuous derivatives. Then, the decomposition of $f$ on a biorthogonal wavelet basis takes the following form:

$$f = \sum_{k \in \mathbb{Z}} \alpha_k \tilde{\phi}_k + \sum_{j \geq 0} \sum_{k \in \mathbb{Z}} \beta_{j,k} \tilde{\psi}_{j,k}, \tag{1.1}$$

where for any $j \geq 0$ and any $k \in \mathbb{Z}$,

$$\alpha_k = \int_{\mathbb{R}} f(x)\phi_k(x)dx, \quad \beta_{j,k} = \int_{\mathbb{R}} f(x)\psi_{j,k}(x)dx.$$

See Section 2.2 for further details. To shorten mathematical expressions, we set

$$\Lambda = \{\lambda = (j,k): \quad j \geq -1, k \in \mathbb{Z}\}$$

and for any $\lambda \in \Lambda$, $\varphi_\lambda = \phi_k$ (respectively $\tilde{\varphi}_\lambda = \tilde{\phi}_k$) if $\lambda = (-1,k)$ and $\varphi_\lambda = \psi_{j,k}$ (respectively $\tilde{\varphi}_\lambda = \tilde{\psi}_{j,k}$) if $\lambda = (j,k)$ with $j \geq 0$. Similarly, $\beta_\lambda = \alpha_k$ if $\lambda = (-1,k)$ and $\beta_\lambda = \beta_{j,k}$ if $\lambda = (j,k)$ with $j \geq 0$. Now, (1.1) can be rewritten as

$$f = \sum_{\lambda \in \Lambda} \beta_\lambda \tilde{\varphi}_\lambda \text{ with } \beta_\lambda = \int \varphi_\lambda(x) f(x) dx. \tag{1.2}$$

In particular, (1.2) holds for the Haar basis where in this case $\tilde{\varphi}_\lambda = \varphi_\lambda$. Now, let us define the thresholding estimate of $f$ by using the properties of Poisson processes. First, we introduce for any $\lambda \in \Lambda$, the natural estimator of $\beta_\lambda$ defined by

$$\hat{\beta}_\lambda = \frac{1}{n} \int \varphi_\lambda(x) dN_x \tag{1.3}$$

that satisfies $\mathbb{E}(\hat{\beta}_\lambda) = \beta_\lambda$. Then, given some parameter $\gamma > 0$, we define the threshold

$$\eta_{\lambda,\gamma} = \sqrt{2\gamma \tilde{V}_{\lambda,n} \log n} + \frac{\gamma \log n}{3n} \|\varphi_\lambda\|_\infty, \tag{1.4}$$

with

$$\tilde{V}_{\lambda,n} = \hat{V}_{\lambda,n} + \sqrt{2\gamma \log n \hat{V}_{\lambda,n} \frac{\|\varphi_\lambda\|_\infty^2}{n^2}} + 3\gamma \log n \frac{\|\varphi_\lambda\|_\infty^2}{n^2}$$

where

$$\hat{V}_{\lambda,n} = \frac{1}{n^2} \int \varphi_\lambda^2(x) dN_x.$$

Note that $\hat{V}_{\lambda,n}$ satisfies $\mathbb{E}(\hat{V}_{\lambda,n}) = V_{\lambda,n}$, where

$$V_{\lambda,n} = \text{Var}(\hat{\beta}_\lambda) = \frac{1}{n} \int \varphi_\lambda^2(x) f(x) dx.$$

Finally given some subset $\Gamma_n$ of $\Lambda$ of the form

$$\Gamma_n = \{\lambda = (j,k) \in \Lambda: \quad j \leq j_0\},$$

where $j_0 = j_0(n)$ is an integer, we set for any $\lambda \in \Lambda$,

$$\tilde{\beta}_\lambda = \hat{\beta}_\lambda 1_{\{|\hat{\beta}_\lambda| \geq \eta_{\lambda,\gamma}\}} 1_{\{\lambda \in \Gamma_n\}}$$

and we set $\tilde{\beta} = (\tilde{\beta}_\lambda)_{\lambda \in \Lambda}$. Finally, the estimator of $f$ is

$$\tilde{f}_{n,\gamma} = \sum_{\lambda \in \Lambda} \tilde{\beta}_\lambda \tilde{\varphi}_\lambda \tag{1.5}$$

and only depends on the choice of $\gamma$ and $j_0$ fixed later. When the Haar basis is used, the estimate is denoted $\tilde{f}_{n,\gamma}^H$ and its wavelet coefficients are denoted $\tilde{\beta}^H = (\tilde{\beta}_\lambda^H)_{\lambda \in \Lambda}$. Thresholding procedures have been introduced by Donoho and Johnstone [17]. The main idea of [17] is that it is sufficient to keep a small amount of the coefficients to have a good estimation of the function $f$. The threshold $\eta_{\lambda,\gamma}$ seems to be defined in a rather complicated manner but is in fact inspired by the universal threshold proposed by [17] in the Gaussian regression framework. The universal threshold of [17] is defined by $\eta_{\lambda,n}^U = \sqrt{2\sigma^2 \log n}$, where $\sigma^2$ (assumed to be known) is the variance of each noisy wavelet coefficient. In our set-up $V_{\lambda,n} = \mathrm{Var}(\hat{\beta}_\lambda)$ depends on $f$, so it is estimated by $\tilde{V}_{\lambda,n}$. Remark that for fixed $\lambda$, when there exists a constant $c_0 > 0$ such that $f(x) \geq c_0$ for $x$ in the support of $\varphi_\lambda$ and if $\|\varphi_\lambda\|_\infty^2 = o_n(n(\log n)^{-1})$, the deterministic term of (1.4) is negligible with respect to the random one and we have asymptotically

$$\eta_{\lambda,\gamma} \approx \sqrt{2\gamma \tilde{V}_{\lambda,n} \log n},$$

which looks like the universal threshold expression if $\gamma$ is close to 1. Actually, the deterministic term of (1.4) allows to consider $\gamma$ close to 1 and to control large deviations terms for high resolution levels. In the same spirit, $V_{\lambda,n}$ is slightly overestimated and we consider $\tilde{V}_{\lambda,n}$ instead of $\hat{V}_{\lambda,n}$ to define the threshold.

The performance of universal thresholding by using the oracle point of view is studied in [17]. In the context of wavelet function estimation by thresholding, the oracle does not tell us the true function, but tells us the coefficients that have to be kept. This "estimator" obtained with the aid of an oracle is not a true estimator, of course, since it depends on $f$. But it represents an ideal for the particular estimation method. The goal of the oracle approach is to derive true estimators which can essentially "mimic" the performance of the "oracle estimator". For Gaussian regression, [17] proved that universal thresholding leads to an estimator that satisfies an oracle inequality: more precisely, the risk of the universal thresholding rule is not larger than the oracle risk up to some logarithmic term which is the price to pay for not having extra information on the locations of the coefficients to keep. So the main question is: does $\tilde{f}_{n,\gamma}$ satisfy a similar oracle inequality? In our framework, it is easy to see that the oracle estimate is $\bar{f} = \sum_{\lambda \in \Gamma_n} \bar{\beta}_\lambda \tilde{\varphi}_\lambda$, where for any $\lambda \in \Gamma_n$, $\bar{\beta}_\lambda = \hat{\beta}_\lambda 1_{\{\beta_\lambda^2 > V_{\lambda,n}\}}$ and we have

$$\mathbb{E}((\bar{\beta}_\lambda - \beta_\lambda)^2) = \min(\beta_\lambda^2, V_{\lambda,n}).$$

By keeping the coefficients $\hat{\beta}_\lambda$ larger than thresholds defined in (1.4), our estimator has a risk that is not larger than the oracle risk, up to a logarithmic term, as stated by the following key result.

**Theorem 1.** *Let us consider a biorthogonal wavelet basis satisfying the properties described in Section 2.2. Let us fix two constants $c \geq 1$ and $c' \in \mathbb{R}$, and let us define for any $n$, $j_0 = j_0(n)$ the integer such that $2^{j_0} \leq n^c (\log n)^{c'} < 2^{j_0+1}$. If $\gamma > c$, then $\tilde{f}_{n,\gamma}$ satisfies the following oracle inequality: for $n$ large enough*

$$\mathbb{E}(\|\tilde{f}_{n,\gamma} - f\|_2^2) \leq C_1 \left[ \sum_{\lambda \in \Gamma_n} \min(\beta_\lambda^2, V_{\lambda,n} \log n) + \sum_{\lambda \notin \Gamma_n} \beta_\lambda^2 \right] + \frac{C_2}{n} \qquad (1.6)$$

*where $C_1$ is a positive constant depending only on $\gamma$, $c$ and the functions that generate the biorthogonal wavelet basis. $C_2$ is also a positive constant depending on $\gamma$, $c$ $c'$, $\|f\|_1$ and the functions that generate the basis.*

Note that Theorem 1 holds with $c = 1$ and $\gamma > 1$. Following the oracle point of view of Donoho and Johnstone, Theorem 1 shows that our procedure is near optimal. The lack of optimality is due to the logarithmic factor. But this term is in some sense unavoidable, as shown later in Theorem 6. Now, let us discuss the near optimality of our procedure from some other perspectives.

## 1.2 Discussion on the assumptions

Previously, we explained why it is crucial to provide theoretical results under very mild assumptions on $f$. Observe that Theorem 1 is established by only assuming that $f$ belongs to $\mathbb{L}_1(\mathbb{R})$ (to ensure that $N_\mathbb{R} < \infty$ almost surely) and $f$ belongs to $\mathbb{L}_2(\mathbb{R})$ (to obtain wavelet decomposition and the study of the performance of $\tilde{f}_{n,\gamma}$ under the $\mathbb{L}_2$-loss). In particular, $f$ can be unbounded and nothing is said about its support which can be unknown or even infinite. The goal of this section is to discuss this last point since, most of the time, estimation is performed by assuming that the intensity has a compact support known by the statistician, usually $[0,1]$. Of course, most of the Poisson data are not generated by an intensity supported by $[0,1]$ and statisticians know this fact but they have in mind a simple preprocessing that can be described as follows. Let us assume that we know a constant $M$ such that the support of $f$ is contained in $[0,M]$. Then, observations are rescaled by dividing each of them by $M$ and new observations (that all depend on $M$) belong to $[0,1]$. An estimator adapted to signals supported by $[0,1]$ can be performed, which leads to a final estimator of $f$ supported by $[0,M]$ by applying the inverse rescaling. Note that such an estimator highly depends on $M$.

Let us go further by describing the situations that may be encountered. If the observations are physical measures given by an instrument that has a limited capacity, then the practitioner usually knows $M$. In this case, if the observations are not concentrated close to 0 but are spread on the whole interval $[0,M]$ in a homogeneous way, then the previous rescaling method performs well. But if one does not have access to $M$ then we are forced in the previous method to estimate it, usually by the largest observation. Then one is forced to face the problem that two different experiments will not lead to estimators with the same support or defined at the same scale and hence it will be hard to compare them. Note also that up to our knowledge, sharp asymptotic properties of such rescaling estimators depending on the largest observation have not been studied. In particular, this method does not seem to be robust if the observations are not compactly supported and if their distribution is heavy-tailed. This situation happens for instance in the financial and geological examples mentioned previously (see [29, 32, 22]) but also in a wide variety of situations (see [12]). In these cases, if observations are rescaled by the largest one, then, methods described at the beginning of the paper provide a very rough estimate of $f$ on small intervals close to 0. However, most of observations may be concentrated close to 0 (for instance for geological data, see [22]) and sharp local estimation at 0 may be of interest. To overcome this problem, statisticians with the help of experts can truncate the data and estimate the intensity on a smaller interval $[0, M_{cut}]$ corresponding to the interval of interest. Then, they face the problem that $M_{cut}$ may be random, subjective, may change from a set of data to another one and may omit values with a potential interest in the future.

So, even if partial solutions exist to overcome issues addressed by the support of $f$, they need a special preprocessing and are not completely justified from a theoretical point of view. We propose a procedure that ignores this preprocessing and which is adapted to non compactly supported Poisson intensities. Our procedure is simple (simpler than the preprocessing described previously) and we prove in the sequel that our method is adaptive minimax with respect to the support which can be bounded or not.

## 1.3 Optimality of $\tilde{f}_{n,\gamma}$ under the minimax approach

To the best of our knowledge, minimax rates for Poisson intensity estimation have not been investigated when the intensity is not compactly supported. But let us mention results established in the following close set-up: the problem of estimating a non-compactly supported density based on the observations of a $n$-sample, which has been partly solved from the minimax point of view. First, let us cite [9] where minimax results for a class of functions depending on a jauge are established or [20] for Sobolev classes. In these papers, the loss function depends on the parameters of the functional class. Similarly, Donoho *et al.* [18] proved the optimality of wavelet linear estimators on Besov spaces $\mathcal{B}_{p,q}^{\alpha}$ when the $\mathbb{L}_p$-risk is considered. First general results where the loss is independent of the functional class have been pointed out by Juditsky and Lambert-Lacroix [25] who investigated minimax rates on the particular class of the Besov spaces $\mathcal{B}_{\infty,\infty}^{\alpha}$ for the $\mathbb{L}_\pi$-risk. When $\pi > 2 + 1/\alpha$, the minimax risk is of the same order up to a logarithmic term as in the equivalent estimation problem on $[0,1]$. However, the behavior of the minimax risk changes dramatically when $\pi \leq 2 + 1/\alpha$, and in this case, it depends on $\pi$. Note that minimax rates for the whole class of Besov spaces $\mathcal{B}_{p,q}^{\alpha}$ ($\alpha > 0$, $1 \leq p, q \leq \infty$) are not derived in [25]. This is the goal of Section 3 under the $\mathbb{L}_2$ risk in the Poisson set-up.

Under mild assumptions on $\gamma$, $\alpha$, $p$, $c$ and $c'$, we prove that the maximal risk of our procedure over balls of $\mathcal{B}_{p,q}^{\alpha} \cap \mathbb{L}_\infty$ is smaller than

$$\left(\frac{\log n}{n}\right)^s \quad \text{with} \quad s = \begin{cases} \frac{2\alpha}{1+2\alpha} & \text{if } 1 \leq p \leq 2 \\ \frac{\alpha}{1+\alpha-\frac{1}{p}} & \text{if } 2 \leq p \leq +\infty. \end{cases}$$

We mention that actually for $p > 2$, it is not necessary to assume that the functions belong to $\mathbb{L}_\infty$ to derive the rate. In addition, we derive the lower bound of the minimax risk for $\mathcal{B}_{p,q}^{\alpha} \cap \mathbb{L}_\infty$ that coincides with the previous upper bound up to the logarithmic term. Let us discuss these results. We note an elbow phenomenon for the rate exponent $s$. When $p \leq 2$, $s$ corresponds to the minimax rate exponent for estimating a compactly supported intensity of a Poisson process. Roughly speaking, it means that it is not harder to estimate non-compactly supported functions than compactly supported functions from the minimax point of view. When $p > 2$, the rate exponent, which depends on $p$, deteriorates when $p$ increases, which means that the support plays a harmful role in this case. An interpretation of this fact and a long discussion of the minimax results are proposed in Section 3.2. Let us just mention that these results are established by using the maxiset approach presented in Section 3.1. We conclude this section by emphasizing that $\tilde{f}_{n,\gamma}$ is rate-optimal, up to the logarithmic term, without knowing the regularity and the support of the underlying signal to be estimated.

## 1.4 Overview of the paper

Section 2 recalls properties of the Poisson process and introduces the biorthogonal wavelet bases used in this paper. Section 3 discusses the properties of our procedure in the minimax and maxiset approaches. Section 4 provides a very general oracle type inequality based on the model selection approach from which Theorem 1 is derived and contains the proofs of the other results.

## 2 Main Tools

### 2.1 Some probabilistic properties of the Poisson process

Let us first recall some basic facts about Poisson processes.

**Definition 1.** *Let $(X, \mathcal{X})$ be a measurable space. Let $N$ be a random countable subset of $X$. $N$ is said to be a Poisson process on $(X, \mathcal{X})$ if*

1. *for any $A \in \mathcal{X}$, the number of points of $N$ lying in $A$ is a random variable, denoted $N_A$, which obeys a Poisson distribution with parameter $\mu(A)$, where $\mu$ is a measure on $X$.*

2. *for any finite family of disjoint sets $A_1, ..., A_n$ of $\mathcal{X}$, $N_{A_1}, ..., N_{A_n}$ are independent random variables.*

The measure $\mu$, called the mean measure of $N$, has no atom (see [27]). In this paper, we assume that $X = \mathbb{R}$, $\mu(\mathbb{R}) < \infty$ and $\mu$ is absolutely continuous with respect with the Lebesgue measure. As explained in Introduction, without loss of generality, we introduce a parameter $n$ and we define the intensity of the process as $f = \frac{d\mu}{n dx}$. We can also mention that a Poisson process $N$ is infinitely divisible, which means that it can be written as follows: for any positive integer $k$:

$$dN = \sum_{i=1}^{k} dN_i \qquad (2.1)$$

where the $N_i$'s are mutually independent Poisson processes on $\mathbb{R}$ with mean measure $\mu/k$. The following proposition (sometimes attributed to Campbell (see [27])) is fundamental and will be used along this paper.

**Proposition 1.** *For any measurable function $g$ and any $z \in \mathbb{R}$, such that $\int e^{zgx} d\mu_x < \infty$ one has,*

$$\mathbb{E}\left[\exp\left(z \int_{\mathbb{R}} g(x) dN_x\right)\right] = \exp\left(\int_{\mathbb{R}} \left(e^{zg(x)} - 1\right) d\mu_x\right).$$

*So,*

$$\mathbb{E}\left(\int_{\mathbb{R}} g(x) dN_x\right) = \int_{\mathbb{R}} g(x) d\mu_x, \quad \text{Var}\left(\int_{\mathbb{R}} g(x) dN_x\right) = \int_{\mathbb{R}} g^2(x) d\mu_x.$$

*If $g$ is bounded, this implies the following exponential inequality. For any $u > 0$,*

$$\mathbb{P}\left(\int_{\mathbb{R}} g(x)(dN_x - d\mu_x) \geq \sqrt{2u \int_{\mathbb{R}} g^2(x) d\mu_x} + \frac{1}{3} \|g\|_\infty u\right) \leq \exp(-u). \qquad (2.2)$$

## 2.2 Biorthogonal wavelet bases and Besov spaces

In this paper, the intensity $f$ to be estimated is assumed to belong to $\mathbb{L}_1 \cap \mathbb{L}_2$. In this case, $f$ can be decomposed on the Haar wavelet basis and this property is used throughout this paper. However, the Haar basis suffers from lack of regularity. To remedy this problem, in particular for deriving minimax properties of $\tilde{f}_{n,\gamma}$ on Besov spaces, we consider a particular class of biorthogonal wavelet bases that are described now. For this purpose, let us set

$$\phi = 1_{[0,1]}.$$

For any $r > 0$, there exist three functions $\psi$, $\tilde{\phi}$ and $\tilde{\psi}$ with the following properties:

1. $\tilde{\phi}$ and $\tilde{\psi}$ are compactly supported,
2. $\tilde{\phi}$ and $\tilde{\psi}$ belong to $C^{r+1}$, where $C^{r+1}$ denotes the Hölder space of order $r + 1$,

3. $\psi$ is compactly supported and is a piecewise constant function,

4. $\psi$ is orthogonal to polynomials of degree no larger than $r$,

5. $\{(\phi_k, \psi_{j,k})_{j\geq 0, k\in\mathbb{Z}}, (\tilde{\phi}_k, \tilde{\psi}_{j,k})_{j\geq 0, k\in\mathbb{Z}}\}$ is a biorthogonal family: for any $j, j' \geq 0$, for any $k, k' \in \mathbb{Z}$,

$$\int_{\mathbb{R}} \psi_{j,k}(x)\tilde{\phi}_{k'}(x)dx = \int_{\mathbb{R}} \phi_k(x)\tilde{\psi}_{j',k'}(x)dx = 0,$$

$$\int_{\mathbb{R}} \phi_k(x)\tilde{\phi}_{k'}(x)dx = 1_{k=k'}, \quad \int_{\mathbb{R}} \psi_{j,k}(x)\tilde{\psi}_{j',k'}(x)dx = 1_{j=j',k=k'},$$

where for any $x \in \mathbb{R}$ and for any $(j, k) \in \mathbb{Z}^2$,

$$\phi_k(x) = \phi(x-k), \quad \psi_{j,k}(x) = 2^{\frac{j}{2}}\psi(2^j x - k)$$

and

$$\tilde{\phi}_k(x) = \tilde{\phi}(x-k), \quad \tilde{\psi}_{j,k}(x) = 2^{\frac{j}{2}}\tilde{\psi}(2^j x - k).$$

This implies the wavelet decomposition (1.1) of $f$. Such biorthogonal wavelet bases have been built by Cohen *et al.* [11] as a special case of spline systems (see also the elegant equivalent construction of Donoho [16] from boxcar functions). The Haar basis can be viewed as a particular biorthogonal wavelet basis, by setting $\tilde{\phi} = \phi$ and $\tilde{\psi} = \psi = 1_{[0,\frac{1}{2}]} - 1_{]\frac{1}{2},1]}$, with $r = 0$ even if Property 2 is not satisfied with such a choice. The Haar basis is an orthonormal basis, which is not true for general biorthogonal wavelet bases. However, we have the frame property: if we denote

$$\Phi = \{\phi, \psi, \tilde{\phi}, \tilde{\psi}\}$$

there exist two constants $c_1(\Phi)$ and $c_2(\Phi)$ only depending on $\Phi$ such that

$$c_1(\Phi)\left(\sum_{k\in\mathbb{Z}}\alpha_k^2 + \sum_{j\geq 0}\sum_{k\in\mathbb{Z}}\beta_{j,k}^2\right) \leq \|f\|_2^2 \leq c_2(\Phi)\left(\sum_{k\in\mathbb{Z}}\alpha_k^2 + \sum_{j\geq 0}\sum_{k\in\mathbb{Z}}\beta_{j,k}^2\right).$$

For instance, when the Haar basis is considered, $c_1(\Phi) = c_2(\Phi) = 1$. In particular, we have

$$c_1(\Phi)\|\tilde{\beta} - \beta\|_{\ell_2}^2 \leq \|\tilde{f}_{n,\gamma} - f\|_2^2 \leq c_2(\Phi)\|\tilde{\beta} - \beta\|_{\ell_2}^2. \tag{2.3}$$

An important feature of such bases is the following: there exists a constant $\mu_\psi > 0$ such that

$$\inf_{x\in[0,1]}|\phi(x)| \geq 1, \quad \inf_{x\in\mathrm{supp}(\psi)}|\psi(x)| \geq \mu_\psi, \tag{2.4}$$

where $\mathrm{supp}(\psi) = \{x \in \mathbb{R}: \psi(x) \neq 0\}$. This property is used throughout the paper.

Now, let us give some properties of Besov spaces that are extensively used in the next section. We recall that Besov spaces, denoted $\mathcal{B}_{p,q}^\alpha$ in the sequel, are defined by using modulus of continuity (see [14] and [21]). We just recall the sequential characterization of Besov spaces by using the biorthogonal wavelet basis (for further details, see [13]).

Let $1 \leq p, q \leq \infty$ and $0 < \alpha < r+1$, the $\mathcal{B}_{p,q}^\alpha$-norm of $f$ is equivalent to the norm

$$\|f\|_{\alpha,p,q} = \begin{cases} \|(\alpha_k)_k\|_{\ell_p} + \left[\sum_{j\geq 0} 2^{jq(\alpha+\frac{1}{2}-\frac{1}{p})}\|(\beta_{j,k})_k\|_{\ell_p}^q\right]^{1/q} & \text{if } q < \infty, \\ \|(\alpha_k)_k\|_{\ell_p} + \sup_{j\geq 0} 2^{j(\alpha+\frac{1}{2}-\frac{1}{p})}\|(\beta_{j,k})_k\|_{\ell_p} & \text{if } q = \infty. \end{cases}$$

We use this norm to define the radius of Besov balls. For any $R > 0$, if $0 < \alpha' \leq \alpha < r+1$, $1 \leq p \leq p' \leq \infty$ and $1 \leq q \leq q' \leq \infty$, we obviously have

$$\mathcal{B}^\alpha_{p,q}(R) \subset \mathcal{B}^\alpha_{p,q'}(R), \quad \mathcal{B}^\alpha_{p,q}(R) \subset \mathcal{B}^{\alpha'}_{p,q}(R).$$

Moreover

$$\mathcal{B}^\alpha_{p,q}(R) \subset \mathcal{B}^{\alpha'}_{p',q}(R) \text{ if } \alpha - \frac{1}{p} \geq \alpha' - \frac{1}{p'}. \tag{2.5}$$

The class of Besov spaces $\mathcal{B}^\alpha_{p,\infty}$ provides a useful tool to classify wavelet decomposed signals with respect to their regularity and sparsity properties (see [24]). Roughly speaking, regularity increases when $\alpha$ increases whereas sparsity increases when $p$ decreases (see Section 3.2).

## 3 Minimax results via the maxiset study

We present in this section the minimax results stated in Introduction. These minimax results are deduced from maxiset results that are first presented. Subsection 3.1 can be omitted on first reading.

### 3.1 The maxiset approach

First, let us describe the maxiset approach which is classical in approximation theory and has been initiated in statistics by Kerkyacharian and Picard [26]. For this purpose, let us assume that we are given $f^*$ an estimation procedure. The maxiset study of $f^*$ consists in deciding the accuracy of $f^*$ by fixing a prescribed rate $\rho^*$ and in pointing out all the functions $f$ such that $f$ can be estimated by the procedure $f^*$ at the target rate $\rho^*$. The maxiset of the procedure $f^*$ for this rate $\rho^*$ is the set of all these functions. More precisely, we restrict our study to the signals belonging to $\mathbb{L}_1 \cap \mathbb{L}_2$ and we set:

**Definition 2.** *Let $\rho^* = (\rho^*_n)_n$ be a decreasing sequence of positive real numbers and let $f^* = (f^*_n)_n$ be an estimation procedure. The maxiset of $f^*$ associated with the rate $\rho^*$ and the $\mathbb{L}_2$-loss is*

$$MS(f^*, \rho^*) = \left\{ f \in \mathbb{L}_1 \cap \mathbb{L}_2 : \sup_n \left[ (\rho^*_n)^{-2} \mathbb{E} \|f^*_n - f\|^2_2 \right] < +\infty \right\},$$

*the ball of radius $R > 0$ of the maxiset is defined by*

$$MS(f^*, \rho^*)(R) = \left\{ f \in \mathbb{L}_1 \cap \mathbb{L}_2 : \sup_n \left[ (\rho^*_n)^{-2} \mathbb{E} \|f^*_n - f\|^2_2 \right] \leq R^2 \right\}.$$

So, the outcome of the maxiset approach is a functional space, which can be viewed as an inversion of the minimax theory where an a priori functional assumption is needed. Obviously, the larger the maxiset, the better the procedure. Maxiset results have been established and extensively discussed in different settings for many classes of estimators and for various rates of convergence. Let us cite for instance [26], [3] and [5] for respectively thresholding rules, Bayes procedures and kernel estimators. More interestingly in our framework, [2] derived maxisets for thresholding rules with data-driven thresholds for density estimation.

The goal of this section is to investigate maxisets for $\tilde{f}_\gamma = (\tilde{f}_{n,\gamma})_n$ and we only focus on rates of the form $\rho_s = (\rho_{n,s})_n$, where $0 < s < \frac{1}{2}$ and for any $n$,

$$\rho_{n,s} = \left(\frac{\log n}{n}\right)^s.$$

So, in the sequel, we investigate for any radius $R > 0$:

$$MS(\tilde{f}_\gamma, \rho_s)(R) = \left\{ f \in \mathbb{L}_1 \cap \mathbb{L}_2 : \quad \sup_n \left[ \left( \frac{\log n}{n} \right)^{-2s} \mathbb{E}\|\tilde{f}_{n,\gamma} - f\|_2^2 \right] \leq R^2 \right\}$$

and to avoid tedious technical aspects related to radius of balls, we use the following notation. If $\mathcal{F}_s$ is a given space

$$MS(\tilde{f}_\gamma, \rho_s) \quad :=: \quad \mathcal{F}_s$$

means in the sequel that for any $R > 0$, there exists $R' > 0$ such that

$$MS(\tilde{f}_\gamma, \rho_s)(R) \cap \mathbb{L}_1(R) \cap \mathbb{L}_2(R) \subset \mathcal{F}_s(R') \cap \mathbb{L}_1(R) \cap \mathbb{L}_2(R)$$

and for any $R' > 0$, there exists $R > 0$ such that

$$\mathcal{F}_s(R') \cap \mathbb{L}_1(R') \cap \mathbb{L}_2(R') \subset MS(\tilde{f}_\gamma, \rho_s)(R) \cap \mathbb{L}_1(R') \cap \mathbb{L}_2(R').$$

To characterize maxisets of $\tilde{f}_\gamma$, we set for any $\lambda \in \Lambda$, $\sigma_\lambda^2 = \int \varphi_\lambda^2(x) f(x) dx$ and we introduce the following spaces.

**Definition 3.** *We define for all $R > 0$ and for all $0 < s < \frac{1}{2}$,*

$$W_s = \left\{ f = \sum_{\lambda \in \Lambda} \beta_\lambda \tilde{\varphi}_\lambda : \quad \sup_{t>0} t^{-4s} \sum_{\lambda \in \Lambda} \beta_\lambda^2 \mathbf{1}_{|\beta_\lambda| \leq \sigma_\lambda t} < \infty \right\},$$

*the ball of radius $R$ associated with $W_s$ is:*

$$W_s(R) = \left\{ f = \sum_{\lambda \in \Lambda} \beta_\lambda \tilde{\varphi}_\lambda : \quad \sup_{t>0} t^{-4s} \sum_{\lambda \in \Lambda} \beta_\lambda^2 \mathbf{1}_{|\beta_\lambda| \leq \sigma_\lambda t} \leq R^{2-4s} \right\},$$

*and for any sequence of spaces $\Gamma = (\Gamma_n)_n$ included in $\Lambda$,*

$$B_{2,\Gamma}^s = \left\{ f = \sum_{\lambda \in \Lambda} \beta_\lambda \tilde{\varphi}_\lambda : \quad \sup_n \left[ \left( \frac{\log n}{n} \right)^{-2s} \sum_{\lambda \notin \Gamma_n} \beta_\lambda^2 \right] < \infty \right\}$$

*and*

$$B_{2,\Gamma}^s(R) = \left\{ f = \sum_{\lambda \in \Lambda} \beta_\lambda \tilde{\varphi}_\lambda : \quad \sup_n \left[ \left( \frac{\log n}{n} \right)^{-2s} \sum_{\lambda \notin \Gamma_n} \beta_\lambda^2 \right] \leq R^2 \right\}.$$

These spaces just depend on the coefficients of the biorthogonal wavelet expansion. In [14], a justification of the form of the radius of $W_s$ and further details are provided. These spaces can be viewed as weak versions of classical Besov spaces, hence they are denoted in the sequel weak Besov spaces. Note that if for all $n$,

$$\Gamma_n = \{\lambda = (j,k) \in \Lambda : \quad j \leq j_0\}$$

with

$$2^{j_0} \leq \left( \frac{n}{\log n} \right)^c < 2^{j_0+1}, \quad c > 0$$

then, $B_{2,\Gamma}^s$ is the classical Besov space $\mathcal{B}_{2,\infty}^{c^{-1}s}$ if the reconstruction wavelets are regular enough. We have the following result.

**Theorem 2.** *Let us fix two constants $c \geq 1$ and $c' \in \mathbb{R}$, and let us define for any $n$, $j_0 = j_0(n)$ the integer such that $2^{j_0} \leq n^c (\log n)^{c'} < 2^{j_0+1}$. Let $\gamma > c$. Then, the procedure defined in (1.5) with the sequence $\Gamma = (\Gamma_n)_n$ such that*

$$\Gamma_n = \{\lambda = (j,k) \in \Lambda : \quad j \leq j_0\}$$

*achieves the following maxiset performance: for all $0 < s < \frac{1}{2}$,*

$$MS(\tilde{f}_\gamma, \rho_s) :=: B^s_{2,\Gamma} \cap W_s.$$

*In particular, if $c' = -c$ and $0 < sc^{-1} < r+1$, where $r$ is the parameter of the biorthogonal basis introduced in Section 2.2,*

$$MS(\tilde{f}_\gamma, \rho_s) :=: \mathcal{B}^{sc^{-1}}_{2,\infty} \cap W_s.$$

The maxiset of $\tilde{f}_\gamma$ is characterized by two spaces: a weak Besov space that is directly connected to the thresholding nature of $\tilde{f}_\gamma$ and the space $B^s_{2,\Gamma}$ that handles the coefficients that are not estimated, which corresponds to the indices $j > j_0$. This maxiset result is similar to the result obtained by Autin [2] in the density estimation setting but our assumptions are less restrictive (see Theorem 5.1 of [2]).

Now, let us point out a family of examples of functions that illustrates the previous result. For this purpose, we only consider the Haar basis that allows to have simple formula for the wavelet coefficients. Let us consider for any $0 < \beta < \frac{1}{2}$, $f_\beta$ such that, for any $x \in \mathbb{R}$,

$$f_\beta(x) = x^{-\beta} 1_{x \in ]0,1]}.$$

The following result points out that if $s$ is small enough, $f_\beta$ belongs to $MS(\tilde{f}_\gamma, \rho_s)$ (so $f_\beta$ can be estimated at the rate $\rho_s$), and in addition $f_\beta \notin \mathbb{L}_\infty$. This result illustrates the fact that the classical assumption $\|f\|_\infty < \infty$ is not necessary to estimate $f$ by our procedure.

**Proposition 2.** *We consider the Haar basis and we set $c' = -c$. For $0 < s < 1/6$, under the assumptions of Theorem 2, if*

$$0 < \beta \leq \frac{1}{2}(1 - 6s),$$

*then for $c$ large enough,*

$$f_\beta \in MS(\tilde{f}^H_\gamma, \rho_s).$$

Let us end this section by explaining the links between maxiset and minimax theories. For this purpose, let $\mathcal{F}$ be a functional space and $\mathcal{F}(R)$ be the ball of radius $R$ associated with $\mathcal{F}$. $\mathcal{F}(R)$ is assumed to be included in a ball of $\mathbb{L}_1 \cap \mathbb{L}_2$. The procedure $\tilde{f}_\gamma$ is said to achieve the rate $\rho_s$ on $\mathcal{F}(R)$ if

$$\sup_n \left[ (\rho_{n,s})^{-2} \sup_{f \in \mathcal{F}(R)} \mathbb{E}\|\tilde{f}_{n,\gamma} - f\|_2^2 \right] < \infty.$$

So, obviously, $\tilde{f}_\gamma$ achieves the rate $\rho_s$ on $\mathcal{F}(R)$ if and only if there exists $R' > 0$ such that

$$\mathcal{F}(R) \subset MS(\tilde{f}_\gamma, \rho_s)(R') \cap \mathbb{L}_1(R') \cap \mathbb{L}_2(R').$$

Using previous results, if $c' = -c$ and if properties of regularity and vanishing moments are satisfied by the wavelet basis, this is satisfied if and only if there exists $R'' > 0$ such that

$$\mathcal{F}(R) \subset \mathcal{B}^{c^{-1}s}_{2,\infty}(R'') \cap W_s(R'') \cap \mathbb{L}_1(R'') \cap \mathbb{L}_2(R''). \tag{3.1}$$

This simple observation will be used to prove some minimax statements of the next section.

### 3.2  Minimax results

To the best of our knowledge, the minimax rate is unknown for $\mathcal{B}_{p,q}^\alpha$ when $p < \infty$. Let us investigate this problem by pointing out the minimax properties of $\tilde{f}_\gamma$ on $\mathcal{B}_{p,q}^\alpha$. For this purpose, we consider the procedure $\tilde{f}_\gamma = (\tilde{f}_{n,\gamma})_n$ defined with

$$\Gamma_n = \{\lambda = (j,k) \in \Lambda : \quad j \leq j_0\}$$

and $j_0 = j_0(n)$ is the integer such that

$$2^{j_0} \leq n^c (\log n)^{-c} < 2^{j_0+1}.$$

The real number $c$ is chosen later. We also set for any $R > 0$,

$$\mathcal{L}_{1,2,\infty}(R) = \{f : \quad \|f\|_1 \leq R, \|f\|_2 \leq R, \|f\|_\infty \leq R\}.$$

In the sequel, minimax results depend on the parameter $r$ of the biorthogonal basis introduced in Section 2.2 to measure the regularity of the reconstruction wavelets $(\tilde{\phi}, \tilde{\psi})$. We first consider the case $p \leq 2$.

**Theorem 3.** *Let $R, R' > 0$, $1 \leq p \leq 2$, $1 \leq q \leq \infty$ and $\alpha \in \mathbb{R}$ such that $\max\left(0, \frac{1}{p} - \frac{1}{2}\right) < \alpha < r+1$. Let $c \geq 1$ large enough such that*

$$\alpha\left(1 - \frac{1}{c(1+2\alpha)}\right) \geq \frac{1}{p} - \frac{1}{2}.$$

*If $\gamma > c$, then for any $n$,*

$$\sup_{f \in \mathcal{B}_{p,q}^\alpha(R) \cap \mathcal{L}_{1,2,\infty}(R')} \mathbb{E}(\|\tilde{f}_{n,\gamma} - f\|_2^2) \leq C(\gamma, c, R, R', \alpha, p, \Phi) \left(\frac{\log n}{n}\right)^{\frac{2\alpha}{2\alpha+1}} \tag{3.2}$$

*where $C(\gamma, c, R, R', \alpha, p, \Phi)$ depends on $R'$, $\gamma$, $c$, on the parameters of the Besov ball and on $\Phi$.*

When $p \leq 2$, the rate of the risk of $\tilde{f}_{n,\gamma}$ corresponds to the minimax rate (up to the logarithmic term) for estimation of a compactly supported intensity of a Poisson process (see [30]), or for estimation of a compactly supported density (see [18]). Roughly speaking, it means that it is not harder to estimate non-compactly supported functions than compactly supported functions from the minimax point of view. In addition, the procedure $\tilde{f}_\gamma$ achieves this classical rate up to a logarithmic term. When $p > 2$ these conclusions do not remain true and we have the following result.

**Theorem 4.** *Let $R, R' > 0$, $2 < p \leq \infty$, $1 \leq q \leq \infty$ and $\alpha \in \mathbb{R}$ such that $0 < \alpha < r+1$. Let $c \geq 1$. If $\gamma > c$, then for any $n$,*

$$\sup_{f \in \mathcal{B}_{p,q}^\alpha(R) \cap \mathbb{L}_1(R') \cap \mathbb{L}_2(R')} \mathbb{E}(\|\tilde{f}_{n,\gamma} - f\|_2^2) \leq C(\gamma, c, R, R', \alpha, p, \Phi) \left(\frac{\log n}{n}\right)^{\frac{\alpha}{\alpha+1-\frac{1}{p}}} \tag{3.3}$$

*where $C(\gamma, c, R, R', \alpha, p, \Phi)$ depends on $R'$, $\gamma$, $c$, on the parameters of the Besov ball and on $\Phi$.*

For $p > 2$, we can note that it is not necessary to assume that signals to be estimated belong to $\mathbb{L}_\infty$ to derive rates of convergence for the risk. Note that when $p = \infty$, the risk is bounded by $\left(\frac{\log n}{n}\right)^{\frac{\alpha}{1+\alpha}}$ up to a constant. In the density estimation setting, this rate was also derived by [25] for their thresholding procedure whose risk was studied on $\mathcal{B}^\alpha_{\infty,\infty}(R)$. Now, combining upper bounds (3.2) and (3.3), for any $R, R' > 0$, $1 \leq p \leq \infty$, $1 \leq q \leq \infty$ and $\alpha \in \mathbb{R}$ such that $\max\left(0, \frac{1}{p} - \frac{1}{2}\right) < \alpha < r+1$, we have:

$$\sup_{f \in \mathcal{B}^\alpha_{p,q}(R) \cap \mathcal{L}_{1,2,\infty}(R')} \mathbb{E}(\|\tilde{f}_{n,\gamma} - f\|_2^2) \leq C(\gamma, c, R, R', \alpha, p, \Phi) \left(\frac{\log n}{n}\right)^{\frac{\alpha}{\alpha + \frac{1}{2} + \left(\frac{1}{2} - \frac{1}{p}\right)_+}}$$

under assumptions of Theorem 3. The following result derives lower bounds of the minimax risk and states that $\tilde{f}_{n,\gamma}$ is rate-optimal up to a logarithmic term.

**Theorem 5.** *Let $R, R' > 0$, $1 \leq p \leq \infty$, $1 \leq q \leq \infty$ and $\alpha \in \mathbb{R}$ such that $\max\left(0, \frac{1}{p} - \frac{1}{2}\right) < \alpha < r + 1$. Then,*

$$\liminf_{n \to +\infty} n^{\frac{\alpha}{\alpha + \frac{1}{2} + \left(\frac{1}{2} - \frac{1}{p}\right)_+}} \inf_{\hat{f}} \sup_{f \in \mathcal{B}^\alpha_{p,q}(R) \cap \mathcal{L}_{1,2,\infty}(R')} \mathbb{E}(\|\hat{f}_n - f\|_2^2) \geq \tilde{C}(\gamma, c, R, R', \alpha, p, \Phi)$$

*where $\tilde{C}(\gamma, c, R, R', \alpha, p, \Phi)$ depends on $R'$, $\gamma$, $c$, on the parameters of the Besov ball and on $\Phi$. Furthermore, let $p^* \geq 1$ and $\alpha^* > 0$ such that*

$$\alpha^*\left(1 - \frac{1}{c(1+2\alpha^*)}\right) \geq \frac{1}{p^*} - \frac{1}{2}. \tag{3.4}$$

*Then, $\tilde{f}_\gamma$ is adaptive minimax up to a logarithmic term on*

$$\left\{\mathcal{B}^\alpha_{p,q}(R) \cap \mathcal{L}_{1,2,\infty}(R'): \quad \alpha^* \leq \alpha < r+1, \; p^* \leq p \leq +\infty, \; 1 \leq q \leq \infty\right\}.$$

Table 1 gathers minimax rates (up to a logarithmic term) obtained for each situation.

|  | $1 \leq p \leq 2$ | $2 \leq p \leq \infty$ |
|---|---|---|
| compact support | $n^{-\frac{2\alpha}{2\alpha+1}}$ | $n^{-\frac{2\alpha}{2\alpha+1}}$ |
| non compact support | $n^{-\frac{2\alpha}{2\alpha+1}}$ | $n^{-\frac{\alpha}{\alpha+1-\frac{1}{p}}}$ |

Table 1: Minimax rates on $\mathcal{B}^\alpha_{p,q} \cap \mathcal{L}_{1,2,\infty}$ (up to a logarithmic term) with $1 \leq p, q \leq \infty$, $\alpha > \max\left(0, \frac{1}{p} - \frac{1}{2}\right)$ under the $\|\cdot\|_2^2$-loss.

Our results show the influence of the support on minimax rates. Note that when restricting on compactly supported signals, when $p > 2$, $\mathcal{B}^\alpha_{p,\infty}(R) \subset \mathcal{B}^\alpha_{2,\infty}(\tilde{R})$ for $\tilde{R}$ large enough and in this case, the rate does not depend on $p$. It is not the case when non-compactly supported signals are considered. Actually, we note an elbow phenomenon at $p = 2$ and the rate deteriorates when $p$ increases. Let us give an interpretation of this observation. Johnstone (1994) showed that when $p < 2$, Besov spaces $\mathcal{B}^\alpha_{p,q}$ model sparse signals where at each level, a very few number of the wavelet coefficients are non-negligible. But these coefficients can be very large. When $p > 2$, $\mathcal{B}_{p,q}^\alpha$-spaces typically model dense signals where the wavelet coefficients are not large but most of them can be non-negligible. This explains why the size of the support plays a role for minimax rates as soon as $p > 2$: when the support is larger, the number of wavelet coefficients to be estimated increases dramatically.

Finally, we note that our procedure achieves the minimax rate, up to a logarithmic term. This logarithmic term is the price we pay for considering thresholding rules. In addition, $\tilde{f}_\gamma$ is near rate-optimal without knowing the regularity and the support of the underlying signal to be estimated.

We end this section by proving that our procedure is adaptive minimax (with the exact exponent of the logarithmic factor) over weak Besov spaces introduced in Section 3.1. For this purpose, we consider signals decomposed on the Haar basis, and we establish the following lower bound with respect to $W_s$. We recall that for any $0 < s < \frac{1}{2}$,

$$\rho_{n,s} = \left(\frac{\log n}{n}\right)^s.$$

**Theorem 6.** *We consider the Haar basis (the spaces $W_s$ and $B_{2,\Gamma}^s$ introduced in Section 3.1 are viewed as sequence spaces). Let*

$$\Gamma_n = \{\lambda = (j,k) \in \Lambda : \quad j \leq j_0\}$$

*with $j_0 = j_0(n)$ the integer such that*

$$2^{j_0} \leq n(\log n)^{-1} < 2^{j_0+1}.$$

*For $0 < s < \frac{1}{2}$ and $R, R', R'' > 0$ such that $R'' \geq 1$ and $R' \geq R^{1-2s} \geq 1$, we have*

$$\liminf_{n \to \infty} \rho_{n,s}^{-2} \inf_{\hat{f}} \sup_{f \in W_s(R) \cap B_{2,\Gamma}^s(R') \cap \mathcal{L}_{1,2,\infty}(R'')} \mathbb{E}(\|\hat{f}_n - f\|_2^2) \geq \tilde{C}(s) R^{2-4s},$$

*where $\tilde{C}(s)$ depends only on $s$ and $\Phi$.*

Using Theorem 2 that provides an upper bound for the risk of our procedure, we immediately deduce the following result.

**Corollary 1.** *The procedure $\tilde{f}_\gamma^H$ defined with*

$$\Gamma_n = \{\lambda = (j,k) \in \Lambda : \quad j \leq j_0\}$$

*with $j_0 = j_0(n)$ the integer such that $2^{j_0} \leq n(\log n)^{-1} < 2^{j_0+1}$ and with $\gamma > 1$ is minimax on $W_s(R) \cap B_{2,\Gamma}^s(R') \cap \mathcal{L}_{1,2,\infty}(R'')$ and is adaptive minimax on*

$$\left\{W_s(R) \cap B_{2,\Gamma}^s(R') \cap \mathcal{L}_{1,2,\infty}(R'') : \quad 0 < s < \frac{1}{2}, 1 \leq R'', \ 1 \leq R \leq R'\right\}.$$

## 4 Proofs via the model selection approach

In this section, we use the model selection approach to provide a very general result with respect to the estimation of a countable family of coefficients. This result is stated in Theorem 7 and is valid for various settings. Applied to the Poisson setting, it allows to establish Theorem 1.

### 4.1 Connections between thresholding and model selection

To describe the model selection approach, let us introduce the following empirical contrast: for any family $\alpha = \{\alpha_\lambda, \lambda \in \Lambda\}$, we set

$$\mathcal{C}_n(\alpha) = -2 \sum_{\lambda \in \Lambda} \alpha_\lambda \hat{\beta}_\lambda + \sum_{\lambda \in \Lambda} \alpha_\lambda^2, \qquad (4.1)$$

which is an unbiased estimator of $\mathcal{C}(\alpha) = \|\beta - \alpha\|_{\ell_2}^2 - \|\beta\|_{\ell_2}^2$. Note that the minimum of $\mathcal{C}$ is achieved for $\alpha = \beta$. Model selection proceeds in two steps: first we consider some family of models $m \subset \Lambda$ and we find $\hat{\beta}(m)$ the minimum of $\mathcal{C}_n$ on each model $m$. Then, we use the data to select a value $\hat{m}$ of $m$ and we take $\hat{\beta}(\hat{m})$ as the final estimator. The first step is immediate in our setting: for any $m \subset \Lambda$,

$$\hat{\beta}(m) = (\hat{\beta}_\lambda 1_{\{\lambda \in m\}})_{\lambda \in \Lambda}$$

and $\mathcal{C}_n(\hat{\beta}(m)) = -\sum_{\lambda \in m} \hat{\beta}_\lambda^2$. Now, the question is : how to choose $\hat{m}$? One could be tempted to choose $m$ as large as possible but this choice would lead to estimates with infinite variance. For this reason, Birgé and Massart [8] proposed to introduce a penalty term associated to each model $m$, denoted $\text{pen}(m)$, and to choose $\hat{m}$ by minimizing

$$\text{Crit}(m) = -\sum_{\lambda \in m} \hat{\beta}_\lambda^2 + \text{pen}(m)$$

over a large class of possible models $m$. For instance, we can fix $\Gamma$ a deterministic subset of $\Lambda$ and consider all the subsets of $\Gamma$. The role of the function $m \to \text{pen}(m)$ is to govern the classical bias-variance tradeoff. Now, if we consider a family of thresholds $(\eta_\lambda)_{\lambda \in \Lambda}$ and if we set for any $m \subset \Gamma$

$$\text{pen}(m) = \sum_{\lambda \in m} \eta_\lambda^2,$$

then the model selection procedure is equivalent to the thresholding rule associated with the family $(\eta_\lambda)_{\lambda \in \Lambda}$:

$$\hat{m} = \{\lambda \in \Gamma : \quad |\hat{\beta}_\lambda| \geq \eta_\lambda\}$$

and

$$\hat{\beta}(\hat{m}) = (\hat{\beta}_\lambda 1_{\{|\hat{\beta}_\lambda| \geq \eta_\lambda\}} 1_{\{\lambda \in \Gamma\}})_{\lambda \in \Lambda} = \tilde{\beta}.$$

Let us note that our method has to be performed for signals with infinite support. So, $\Gamma$ may be infinite, which is not usual in the literature. The following theorem is self-contained; we do not use the Poisson setting and we do not make any assumption on the distribution of $\hat{\beta}_\lambda$ or on the form of the threshold $\eta_\lambda$. So, Theorem 7 can be used for other settings and this is the main reason for the following very abstract formulation.

**Theorem 7.** *To estimate a countable family $\beta = (\beta_\lambda)_{\lambda \in \Lambda}$, such that $\|\beta\|_{\ell_2} < \infty$, we assume that a family of coefficient estimators $(\hat{\beta}_\lambda)_{\lambda \in \Gamma}$, where $\Gamma$ is a known deterministic subset of $\Lambda$, and a family of possibly random thresholds $(\eta_\lambda)_{\lambda \in \Gamma}$ are available and we consider the thresholding rule $\tilde{\beta} = (\hat{\beta}_\lambda 1_{|\hat{\beta}_\lambda| \geq \eta_\lambda} 1_{\lambda \in \Gamma})_{\lambda \in \Lambda}$. Let $\varepsilon > 0$ be fixed. Assume that there exist a deterministic family $(F_\lambda)_{\lambda \in \Gamma}$ and three constants $\kappa \in [0, 1[$, $\omega \in [0, 1]$ and $\mu > 0$ (that may depend on $\varepsilon$ but not on $\lambda$) with the following properties.*

*(A1) For all $\lambda$ in $\Gamma$,*

$$\mathbb{P}(|\hat{\beta}_\lambda - \beta_\lambda| > \kappa \eta_\lambda) \leq \omega.$$

*(A2) There exist $1 < p, q < \infty$ with $\frac{1}{p} + \frac{1}{q} = 1$ and a constant $R > 0$ such that for all $\lambda$ in $\Gamma$,*

$$\left(\mathbb{E}(|\hat{\beta}_\lambda - \beta_\lambda|^{2p})\right)^{\frac{1}{p}} \leq R \max(F_\lambda, F_\lambda^{\frac{1}{p}} \varepsilon^{\frac{1}{q}}).$$

*(A3) There exists a constant $\theta$ such that for all $\lambda$ in $\Gamma$ such that $F_\lambda < \theta \varepsilon$*

$$\mathbb{P}(|\hat{\beta}_\lambda - \beta_\lambda| > \kappa \eta_\lambda, |\hat{\beta}_\lambda| > \eta_\lambda) \leq F_\lambda \mu.$$

*Then the estimator $\tilde{\beta}$ satisfies*

$$\frac{1-\kappa^2}{1+\kappa^2} \mathbb{E}\|\tilde{\beta} - \beta\|_{\ell_2}^2 \leq \mathbb{E} \inf_{m \subset \Gamma} \left\{ \frac{1+\kappa^2}{1-\kappa^2} \sum_{\lambda \notin m} \beta_\lambda^2 + \frac{1-\kappa^2}{\kappa^2} \sum_{\lambda \in m} (\hat{\beta}_\lambda - \beta_\lambda)^2 + \sum_{\lambda \in m} \eta_\lambda^2 \right\} + LD \sum_{\lambda \in \Gamma} F_\lambda$$

*with*

$$LD = \frac{R}{\kappa^2} \left( \left(1 + \theta^{-1/q}\right) \omega^{1/q} + (1 + \theta^{1/q}) \varepsilon^{1/q} \mu^{1/q} \right).$$

Observe that this result makes sense only when $\sum_{\lambda \in \Gamma} F_\lambda < \infty$ and in this case, if $LD$ (which stands for large deviation inequalities) is small enough, the main term of the right hand side is given by the first term.

Now, let us briefly comment the assumptions of this theorem. The concentration inequality of Assumption (A1) controls the deviation of $|\hat{\beta}_\lambda - \beta_\lambda|$ with respect to 0. The family $(F_\lambda)_{\lambda \in \Gamma}$ is introduced for Assumptions (A2) and (A3). Assumption (A2) provides upper bounds for the moments of $\hat{\beta}_\lambda$ and looks like a Rosenthal inequality if $F_\lambda$ can be related to the variance of $\hat{\beta}_\lambda$. Actually, compactly supported signals can be well estimated by thresholding if sharp concentration and Rosenthal inequalities are satisfied (see Theorem 3 of [18] and Theorem 3.1. of [26]). In our set-up where the support of $f$ can be infinite, these basic tools are not sufficient and Assumption (A3) is introduced to ensure that with high probability, when $F_\lambda$ is small, then either $\beta_\lambda$ is estimated by 0, or $|\hat{\beta}_\lambda - \beta_\lambda|$ is small. Remark 1 in Section 4.2 provides additional technical reasons for the introduction of Assumption (A3) when the support of the signal is infinite. Finally, the condition $\sum_{\lambda \in \Gamma} F_\lambda < \infty$ shows that the variations of $(\hat{\beta}_\lambda)_{\lambda \in \Gamma}$ around $(\beta_\lambda)_{\lambda \in \Gamma}$, as pointed out by Assumptions (A2) and (A3), have to be controlled in a global way.

This theorem applied in the Poisson set-up with $\Gamma = \Gamma_n$ and $\eta_\lambda = \eta_{\lambda,\gamma}$ implies Theorem 1. In particular the family $(F_\lambda)_\lambda$ is given by $F_\lambda = \int_{\text{supp}(\varphi_\lambda)} f(x) dx$, which is related to the variance of $\hat{\beta}_\lambda$ (see (4.5)).

Using (2.3), without loss of generality, Theorems 1, 2, 3, 4 and 5 are established by using the $\ell_2$-norm of coefficients instead of the functional $\mathbb{L}_2$-loss. In the following proofs, the values of the constants $C_1, C_2, K_1, K_2, \theta, \ldots$ may change from one proof to another one. Finally, recall that we have set for any $\lambda \in \Lambda$,

$$\sigma_\lambda^2 = \int \varphi_\lambda^2(x) f(x) dx.$$

## 4.2 Proof of Theorem 7

We use the model selection approach. By definition of $\hat{m}$ one has for any $m \subset \Gamma$,

$$\mathcal{C}_n(\tilde{\beta}) + \text{pen}(\hat{m}) \leq \mathcal{C}_n(\hat{\beta}(m)) + \text{pen}(m).$$

For any family $\alpha = (\alpha_\lambda)_{\lambda \in \Lambda}$, we set

$$\nu(\alpha) = \sum_{\lambda \in \Lambda} \alpha_\lambda (\hat{\beta}_\lambda - \beta_\lambda).$$

Then, using (4.1),

$$\mathcal{C}_n(\alpha) = \|\beta - \alpha\|_{\ell_2}^2 - \|\beta\|_{\ell_2}^2 - 2\nu(\alpha).$$

So,

$$\begin{aligned}
\|\tilde{\beta} - \beta\|_{\ell_2}^2 &\leq \|\hat{\beta}(m) - \beta\|_{\ell_2}^2 + 2\nu(\tilde{\beta} - \hat{\beta}(m)) + \text{pen}(m) - \text{pen}(\hat{m}) \\
&\leq \|\hat{\beta}(m) - \beta\|_{\ell_2}^2 + 2\nu(\tilde{\beta} - \beta(m)) - 2\nu(\hat{\beta}(m) - \beta(m)) + \text{pen}(m) - \text{pen}(\hat{m}),
\end{aligned}$$

where $\beta(m) = \mathbb{E}(\hat{\beta}(m))$ is the projection of $\beta$ on the space of the vectors $\alpha = (\alpha_\lambda)_{\lambda \in \Lambda}$ such that $\alpha_\lambda = 0$ when $\lambda \notin m$ for the $\ell_2$-norm. But,

$$\|\hat{\beta}(m) - \beta\|_{\ell_2}^2 = \|\hat{\beta}(m) - \beta(m)\|_{\ell_2}^2 + \|\beta(m) - \beta\|_{\ell_2}^2 = \nu(\hat{\beta}(m) - \beta(m)) + \|\beta(m) - \beta\|_{\ell_2}^2$$

and

$$\begin{aligned}
2\nu(\tilde{\beta} - \beta(m)) &\leq 2\|\tilde{\beta} - \beta(m)\|_{\ell_2} \chi(m \cup \hat{m}) \\
&\leq 2\|\tilde{\beta} - \beta\|_{\ell_2} \chi(m \cup \hat{m}) + 2\|\beta - \beta(m)\|_{\ell_2} \chi(m \cup \hat{m}) \\
&\leq \frac{2\kappa^2}{1+\kappa^2} \|\tilde{\beta} - \beta\|_{\ell_2}^2 + \frac{2\kappa^2}{1-\kappa^2} \|\beta(m) - \beta\|_{\ell_2}^2 + \frac{1}{\kappa^2} \chi^2(m \cup \hat{m}),
\end{aligned}$$

where we have set for any $m \subset \Gamma$,

$$\chi(m) = \|\hat{\beta}(m) - \beta(m)\|_{\ell_2} = \sqrt{\sum_{\lambda \in m} (\hat{\beta}_\lambda - \beta_\lambda)^2} = \sqrt{\nu(\hat{\beta}(m) - \beta(m))}$$

and we have used twice the inequality $2ab \leq \rho a^2 + \rho^{-1} b^2$ with $\rho = 2\kappa^2(1+\kappa^2)^{-1}$ and $\rho = 2\kappa^2(1-\kappa^2)^{-1}$. Finally,

$$\begin{aligned}
\frac{1-\kappa^2}{1+\kappa^2} \|\tilde{\beta} - \beta\|_{\ell_2} &\leq -\|\hat{\beta}(m) - \beta(m)\|_{\ell_2}^2 + \frac{1+\kappa^2}{1-\kappa^2} \|\beta(m) - \beta\|_{\ell_2}^2 + \frac{1}{\kappa^2} \chi^2(m \cup \hat{m}) + \text{pen}(m) - \text{pen}(\hat{m}) \\
&\leq \frac{1+\kappa^2}{1-\kappa^2} \|\beta(m) - \beta\|_{\ell_2}^2 + \left(\frac{1}{\kappa^2} - 1\right) \|\hat{\beta}(m) - \beta(m)\|_{\ell_2}^2 + \text{pen}(m) + \mathcal{A},
\end{aligned}$$

where

$$\mathcal{A} = \frac{1}{\kappa^2} \chi^2(\hat{m}) - \text{pen}(\hat{m}) = \sum_{\lambda \in \Gamma} \left(\frac{1}{\kappa^2} (\hat{\beta}_\lambda - \beta_\lambda)^2 - \eta_\lambda^2\right) 1_{|\hat{\beta}_\lambda| > \eta_\lambda}.$$

Now, we introduce

$$A_1 = \sum_{\lambda \in \Gamma} \mathbb{E}\left(\frac{1}{\kappa^2} (\hat{\beta}_\lambda - \beta_\lambda)^2 1_{|\hat{\beta}_\lambda - \beta_\lambda| > \kappa \eta_\lambda}\right) 1_{F_\lambda \geq \theta \varepsilon}$$

and

$$A_2 = \sum_{\lambda \in \Gamma} \mathbb{E}\left(\frac{1}{\kappa^2} (\hat{\beta}_\lambda - \beta_\lambda)^2 1_{|\hat{\beta}_\lambda - \beta_\lambda| > \kappa \eta_\lambda} 1_{|\hat{\beta}_\lambda| > \eta_\lambda}\right) 1_{F_\lambda < \theta \varepsilon}.$$

Therefore,

$$\mathbb{E}[\mathcal{A}] \leq A_1 + A_2.$$

By using the Hölder inequality,

$$\begin{aligned}
A_1 &\leq \frac{1}{\kappa^2}\sum_{\lambda\in\Gamma}\left(\mathbb{E}(|\hat{\beta}_\lambda-\beta_\lambda|^{2p})\right)^{\frac{1}{p}}\left(\mathbb{P}(|\hat{\beta}_\lambda-\beta_\lambda|>\kappa\eta_\lambda)\right)^{\frac{1}{q}}1_{F_\lambda\geq\theta\varepsilon}\\
&\leq \frac{R\omega^{\frac{1}{q}}}{\kappa^2}\sum_{\lambda\in\Gamma}\max(F_\lambda,F_\lambda^{\frac{1}{p}}\varepsilon^{\frac{1}{q}})1_{F_\lambda\geq\theta\varepsilon}\\
&\leq \frac{R\omega^{\frac{1}{q}}}{\kappa^2}\left(\sum_{\lambda\in\Gamma}F_\lambda+\varepsilon^{\frac{1}{q}}\sum_{\lambda\in\Gamma}F_\lambda^{\frac{1}{p}}\left(\frac{F_\lambda}{\theta\varepsilon}\right)^{\frac{1}{q}}\right)\\
&\leq \frac{R\omega^{\frac{1}{q}}}{\kappa^2}\left(1+\theta^{-\frac{1}{q}}\right)\sum_{\lambda\in\Gamma}F_\lambda
\end{aligned}$$

and

$$\begin{aligned}
A_2 &\leq \frac{1}{\kappa^2}\sum_{\lambda\in\Gamma}\left(\mathbb{E}(|\hat{\beta}_\lambda-\beta_\lambda|^{2p})\right)^{\frac{1}{p}}\left(\mathbb{P}(|\hat{\beta}_\lambda-\beta_\lambda|>\kappa\eta_\lambda,|\hat{\beta}_\lambda|>\eta_\lambda)\right)^{\frac{1}{q}}1_{F_\lambda<\theta\varepsilon}\\
&\leq \frac{R}{\kappa^2}\sum_{\lambda\in\Gamma}\max(F_\lambda,F_\lambda^{\frac{1}{p}}\varepsilon^{\frac{1}{q}})F_\lambda^{\frac{1}{q}}\mu^{\frac{1}{q}}1_{F_\lambda<\theta\varepsilon}\\
&\leq \frac{R}{\kappa^2}\left(\sum_{\lambda\in\Gamma}F_\lambda^{1+\frac{1}{q}}\mu^{\frac{1}{q}}\left(\frac{\theta\varepsilon}{F_\lambda}\right)^{\frac{1}{q}}+\mu^{\frac{1}{q}}\varepsilon^{\frac{1}{q}}\sum_{\lambda\in\Gamma}F_\lambda\right)\\
&\leq \frac{R}{\kappa^2}\mu^{\frac{1}{q}}(1+\theta^{\frac{1}{q}})\varepsilon^{\frac{1}{q}}\sum_{\lambda\in\Gamma}F_\lambda.
\end{aligned}$$

So,

$$\mathbb{E}(\mathcal{A})\leq LD\sum_{\lambda\in\Gamma}F_\lambda,$$

which proves Theorem 7.

**Remark 1.** *When compactly supported signals are considered, it is natural to take $\Gamma$ satisfying $card(\Gamma)<\infty$ and in this case, the upper bound of $\mathbb{E}(\mathcal{A})$ takes the simpler form:*

$$\begin{aligned}
\mathbb{E}(\mathcal{A}) &\leq \frac{1}{\kappa^2}\sum_{\lambda\in\Gamma}\left(\mathbb{E}(|\hat{\beta}_\lambda-\beta_\lambda|^{2p})\right)^{\frac{1}{p}}\left(\mathbb{P}(|\hat{\beta}_\lambda-\beta_\lambda|>\kappa\eta_\lambda)\right)^{\frac{1}{q}}\\
&\leq \frac{1}{\kappa^2}card(\Gamma)\max_{\lambda\in\Gamma}\left(\mathbb{E}(|\hat{\beta}_\lambda-\beta_\lambda|^{2p})\right)^{\frac{1}{p}}w^{\frac{1}{q}}.
\end{aligned}$$

*Even under a rough control of $\max_{\lambda\in\Gamma}\mathbb{E}(|\hat{\beta}_\lambda-\beta_\lambda|^{2p})$, the term $\mathbb{E}(\mathcal{A})$ is negligible with respect to the main term as soon as $w$ is small enough, which occurs if the threshold is large enough. In particular, when restricting our attention to compactly supported signals, Assumption (A3) is useless.*

### 4.3 Proof of Theorem 1

To prove Theorem 1, we use Theorem 7 with $\hat{\beta}_\lambda$ defined in (1.3), $\eta_\lambda=\eta_{\lambda,\gamma}$ defined in (1.4) and

$$\Gamma=\Gamma_n=\{\lambda=(j,k)\in\Lambda:\quad -1\leq j\leq j_0\}\text{ with }2^{j_0}\leq n^c(\log n)^{c'}<2^{j_0+1}.$$

We set
$$F_\lambda = \int_{\mathrm{supp}(\varphi_\lambda)} f(x)dx,$$
so we have:
$$\sum_{\lambda \in \Gamma} F_\lambda = \sum_{-1 \leq j \leq j_0} \sum_k \int_{x \in \mathrm{supp}(\varphi_{j,k})} f(x)dx \leq \int f(x)dx \sum_{-1 \leq j \leq j_0} \sum_k 1_{x \in \mathrm{supp}(\varphi_{j,k})} \leq (j_0+2)m_\varphi \|f\|_1, \quad (4.2)$$
where $m_\varphi$ is a finite constant depending only on the compactly supported functions $\phi$ and $\psi$. Finally, $\sum_{\lambda \in \Gamma} F_\lambda$ is bounded by $\log(n)$ up to a constant that only depends on $\|f\|_1$, $c$, $c'$ and the functions $\phi$ and $\psi$. Now, we give a fundamental lemma to derive Assumption (A1) of Theorem 7.

**Lemma 1.** *For any $u > 0$*
$$\mathbb{P}\left(|\hat{\beta}_\lambda - \beta_\lambda| \geq \sqrt{2uV_{\lambda,n}} + \frac{\|\varphi_\lambda\|_\infty u}{3n}\right) \leq 2e^{-u}. \quad (4.3)$$

*Moreover, for any $u > 0$*
$$\mathbb{P}\left(V_{\lambda,n} \geq \breve{V}_{\lambda,n}(u)\right) \leq e^{-u},$$
*where*
$$\breve{V}_{\lambda,n}(u) = \hat{V}_{\lambda,n} + \sqrt{2\hat{V}_{\lambda,n}\frac{\|\varphi_\lambda\|_\infty^2}{n^2}u} + 3\frac{\|\varphi_\lambda\|_\infty^2}{n^2}u.$$

**Proof.** Equation (4.3) comes easily from (2.2) applied with $g = \varphi_\lambda/n$. The same inequality applied with $g = -\varphi_\lambda^2/n^2$ gives:
$$\mathbb{P}\left(V_{\lambda,n} \geq \hat{V}_{\lambda,n} + \sqrt{2u\int_\mathbb{R} \frac{\varphi_\lambda^4(x)}{n^4}nf(x)dx} + \frac{\|\varphi_\lambda\|_\infty^2}{3n^2}u\right) \leq e^{-u}.$$

We observe that
$$\int_\mathbb{R} \frac{\varphi_\lambda^4(x)}{n^4}nf(x)dx \leq \frac{\|\varphi_\lambda\|_\infty^2}{n^2}V_{\lambda,n}.$$

So, if we set $a = u\frac{\|\varphi_\lambda\|_\infty^2}{n^2}$, then
$$\mathbb{P}(V_{\lambda,n} - \sqrt{2V_{\lambda,n}a} - a/3 \geq \hat{V}_{\lambda,n}) \leq e^{-u}.$$

We obtain
$$\mathbb{P}(\sqrt{V_{\lambda,n}} \geq \mathcal{P}^{-1}(\hat{V}_{\lambda,n})) \leq e^{-u}$$
where $\mathcal{P}^{-1}(\hat{V}_{\lambda,n})$ is the positive solution of
$$(\mathcal{P}^{-1}(\hat{V}_{\lambda,n}))^2 - \sqrt{2a}\mathcal{P}^{-1}(\hat{V}_{\lambda,n}) - (a/3 + \hat{V}_{\lambda,n}) = 0.$$

To conclude, it remains to observe that
$$\breve{V}_{\lambda,n}(u) \geq (\mathcal{P}^{-1}(\hat{V}_{\lambda,n}))^2 = \left(\sqrt{\hat{V}_{\lambda,n} + 5a/6} + \sqrt{a/2}\right)^2.$$

∎

Let $\kappa < 1$. Combining these inequalities with $\tilde{V}_{\lambda,n} = \check{V}_{\lambda,n}(\gamma \log n)$ yields

$$
\begin{aligned}
\mathbb{P}(|\hat{\beta}_\lambda - \beta_\lambda| > \kappa \eta_{\lambda,\gamma}) &\leq \mathbb{P}\left(|\hat{\beta}_\lambda - \beta_\lambda| \geq \sqrt{2\kappa^2 \gamma \log n \tilde{V}_{\lambda,n}} + \frac{\kappa \gamma \log n \|\varphi_\lambda\|_\infty}{3n}\right) \\
&\leq \mathbb{P}\left(|\hat{\beta}_\lambda - \beta_\lambda| \geq \sqrt{2\kappa^2 \gamma \log n \tilde{V}_{\lambda,n}} + \frac{\kappa \gamma \log n \|\varphi_\lambda\|_\infty}{3n}, V_{\lambda,n} \geq \tilde{V}_{\lambda,n}\right) \\
&\quad + \mathbb{P}\left(|\hat{\beta}_\lambda - \beta_\lambda| \geq \sqrt{2\kappa^2 \gamma \log n \tilde{V}_{\lambda,n}} + \frac{\kappa \gamma \log n \|\varphi_\lambda\|_\infty}{3n}, V_{\lambda,n} < \tilde{V}_{\lambda,n}\right) \\
&\leq \mathbb{P}(V_{\lambda,n} \geq \tilde{V}_{\lambda,n}) + \mathbb{P}\left(|\hat{\beta}_\lambda - \beta_\lambda| \geq \sqrt{2\kappa^2 \gamma \log n V_{\lambda,n}} + \frac{\kappa \gamma \log n \|\varphi_\lambda\|_\infty}{3n}\right) \\
&\leq n^{-\gamma} + 2n^{-\kappa^2 \gamma} \\
&\leq 3n^{-\kappa^2 \gamma}.
\end{aligned}
$$

So, for any value of $\kappa \in [0, 1[$, Assumption (A1) is true with $\eta_\lambda = \eta_{\lambda,\gamma}$ if we take $\omega = 3n^{-\kappa^2 \gamma}$. To verify the Rosenthal type inequality (A2) of Theorem 7, we prove the following lemma.

**Lemma 2.** *For any $p \geq 2$, there exists an absolute constant $C$ such that*

$$
\mathbb{E}(|\hat{\beta}_\lambda - \beta_\lambda|^{2p}) \leq C^p p^{2p} \left(V_{\lambda,n}^p + \left[\frac{\|\varphi_\lambda\|_\infty}{n}\right]^{2p-2} V_{\lambda,n}\right).
$$

**Proof.** We apply (2.1). Hence,

$$
\hat{\beta}_\lambda - \beta_\lambda = \sum_{i=1}^k \int \frac{\varphi_\lambda(x)}{n} \left(dN_x^i - nk^{-1} f(x) dx\right) = \sum_{i=1}^k Y_i
$$

where for any $i$,

$$
Y_i = \int \frac{\varphi_\lambda(x)}{n} \left(dN_x^i - nk^{-1} f(x) dx\right).
$$

So the $Y_i$'s are i.i.d. centered variables, each of them has a moment of order $2p$. For any $i$, we apply the Rosenthal inequality (see Theorem 2.5 of [23]) to the positive and negative parts of $Y_i$. This easily implies that

$$
\mathbb{E}\left(\left|\sum_{i=1}^k Y_i\right|^{2p}\right) \leq \left(\frac{16p}{\log(2p)}\right)^{2p} \max\left(\left(\mathbb{E} \sum_{i=1}^k Y_i^2\right)^p, \left(\mathbb{E} \sum_{i=1}^k |Y_i|^{2p}\right)\right).
$$

It remains to bound the upper limit of $\mathbb{E}(\sum_{i=1}^k |Y_i|^\ell)$ for all $\ell \in \{2p, 2\} \geq 2$ when $k \to \infty$. Let us introduce

$$
\Omega_k = \{\forall\, i \in \{1, \ldots, k\}, N_\mathbb{R}^i \leq 1\}.
$$

Then, it is easy to see that $\mathbb{P}(\Omega_k^c) \leq k^{-1}(n\|f\|_1)^2$ (see e.g., (4.6) below).

On $\Omega_k$, $|Y_i|^\ell = O_k(k^{-\ell})$ if $\int \frac{\varphi_\lambda(x)}{n} dN_x^i = 0$ and $|Y_i|^\ell = \left[\frac{|\varphi_\lambda(T)|}{n}\right]^\ell + O_k\left(k^{-1}\left[\frac{|\varphi_\lambda(T)|}{n}\right]^{\ell-1}\right)$ if

$\int \frac{\varphi_\lambda(x)}{n} dN_x^i = \frac{\varphi_\lambda(T)}{n}$ where $T$ is the point of the process $N^i$. Consequently,

$$\mathbb{E}\sum_{i=1}^k |Y_i|^\ell \leq \mathbb{E}\left(1_{\Omega_k}\left(\sum_{T\in N}\left[\left[\frac{|\varphi_\lambda(T)|}{n}\right]^\ell + O_k\left(k^{-1}\left[\frac{|\varphi_\lambda(T)|}{n}\right]^{\ell-1}\right)\right] + kO_k(k^{-\ell})\right)\right)$$
$$+ \sqrt{\mathbb{P}(\Omega_k^c)}\sqrt{\mathbb{E}\left[\left(\sum_{i=1}^k |Y_i|^\ell\right)^2\right]}. \quad (4.4)$$

But we have

$$\sum_{i=1}^k |Y_i|^\ell \leq 2^{\ell-1}\left(\sum_{i=1}^k\left[\left[\frac{\|\varphi_\lambda\|_\infty}{n}\right]^\ell (N_\mathbb{R}^i)^\ell + \left(k^{-1}\int |\varphi_\lambda(x)|f(x)dx\right)^\ell\right]\right)$$
$$\leq 2^{\ell-1}\left(\left[\frac{\|\varphi_\lambda\|_\infty}{n}\right]^\ell N_\mathbb{R}^\ell + k\left(k^{-1}\int |\varphi_\lambda(x)|f(x)dx\right)^\ell\right).$$

So, when $k \to +\infty$, the last term in (4.4) converges to 0 since a Poisson variable has moments of every order and

$$\limsup_{k\to\infty} \mathbb{E}\sum_{i=1}^k |Y_i|^\ell \leq \mathbb{E}\left(\int \left[\frac{|\varphi_\lambda(x)|}{n}\right]^\ell dN_x\right) \leq \left[\frac{\|\varphi_\lambda\|_\infty}{n}\right]^{\ell-2} V_{\lambda,n},$$

which concludes the proof. ∎

Now,

$$V_{\lambda,n} = \frac{1}{n}\int \varphi_\lambda^2(x)f(x)dx \leq \frac{\|\varphi_\lambda\|_\infty^2 F_\lambda}{n} \quad (4.5)$$

and Assumption (A2) is satisfied with $\varepsilon = \frac{1}{n}$ and

$$R = \frac{2Cp^2 2^{j_0}\max(\|\phi\|_\infty^2; \|\psi\|_\infty^2)}{n}$$

since $\|\varphi_\lambda\|_\infty^2 \leq 2^{j_0}\max(\|\phi\|_\infty^2; \|\psi\|_\infty^2)$ and

$$\left(\mathbb{E}(|\hat{\beta}_\lambda - \beta_\lambda|^{2p})\right)^{\frac{1}{p}} \leq Cp^2\left(\frac{\|\varphi_\lambda\|_\infty^2 F_\lambda}{n} + \|\varphi_\lambda\|_\infty^2 F_\lambda^{\frac{1}{p}} n^{\frac{1}{p}-2}\right) \leq \frac{Cp^2\|\varphi_\lambda\|_\infty^2}{n}\left(F_\lambda + F_\lambda^{\frac{1}{p}} n^{-\frac{1}{q}}\right).$$

Finally, Assumption (A3) comes from the following lemma.

**Lemma 3.** *We set*

$$N_\lambda = \int_{supp(\varphi_\lambda)} dN \quad and \quad C' = (\sqrt{6}+1/3)\gamma \geq \sqrt{6}+1/3.$$

*There exists an absolute constant $0 < \theta' < 1$ such that if $nF_\lambda \leq \theta' C'\log n$ and $(1-\theta')(\sqrt{6}+1/3)\log n \geq 2$ then,*

$$\mathbb{P}(N_\lambda - nF_\lambda \geq (1-\theta')C'\log n) \leq F_\lambda n^{-\gamma}.$$

**Remark 2.** *We can take $\theta' = 0.01$ and in this case, the result is true as soon as $n \geq 3$*

**Proof.** One takes $\theta' \in [0, 1]$ (for instance $\theta' = 0.01$) such that

$$\frac{3(1-\theta')^2}{2(2\theta'+1)}(\sqrt{6}+1/3) \geq 4.$$

We use Equation (5.2) of [30] to obtain

$$\mathbb{P}(N_\lambda - nF_\lambda \geq (1-\theta')C'\log n) \leq \exp\left(-\frac{((1-\theta')C'\log n)^2}{2(nF_\lambda + (1-\theta')C'\log n/3)}\right) \leq n^{-\frac{3(1-\theta')^2}{2(2\theta'+1)}C'}.$$

If $nF_\lambda \geq n^{-\gamma-1}$, since $\frac{3(1-\theta')^2}{2(2\theta'+1)}C' \geq 2\gamma + 2$, the result is true. If $nF_\lambda \leq n^{-\gamma-1}$,

$$\mathbb{P}(N_\lambda - nF_\lambda \geq (1-\theta')C'\log n) \leq \mathbb{P}(N_\lambda > (1-\theta')C'\log n) \leq \mathbb{P}(N_\lambda \geq 2) \leq \sum_{k \geq 2} \frac{(nF_\lambda)^k}{k!}e^{-nF_\lambda} \leq (nF_\lambda)^2 \tag{4.6}$$

and the result is true. ∎

Now, observe that if $|\hat{\beta}_\lambda| > \eta_{\lambda,\gamma}$ then

$$N_\lambda \geq C'\log n.$$

Indeed, $|\hat{\beta}_\lambda| > \eta_{\lambda,\gamma}$ implies

$$\frac{C'\log n}{n}\|\varphi_\lambda\|_\infty \leq |\hat{\beta}_\lambda| \leq \frac{\|\varphi_\lambda\|_\infty N_\lambda}{n}.$$

So if $n$ satisfies $(1-\theta')(\sqrt{6}+1/3)\log n \geq 2$, we set $\theta = \theta'C'\log(n)$ and $\mu = n^{-\gamma}$. In this case, Assumption (A3) is fulfilled since if $nF_\lambda \leq \theta'C'\log n$

$$\mathbb{P}(|\hat{\beta}_\lambda - \beta_\lambda| > \kappa\eta_\lambda, |\hat{\beta}_\lambda| > \eta_\lambda) \leq \mathbb{P}(N_\lambda - nF_\lambda \geq (1-\theta')C'\log n) \leq F_\lambda n^{-\gamma}.$$

Finally, if $n$ satisfies $(1-\theta')(\sqrt{6}+1/3)\log n \geq 2$, Theorem 7 applies:

$$\frac{1-\kappa^2}{1+\kappa^2}\mathbb{E}\|\tilde{\beta} - \beta\|_{\ell_2}^2 \leq \inf_{m \subset \Gamma}\left\{\frac{1+\kappa^2}{1-\kappa^2}\sum_{\lambda \notin m}\beta_\lambda^2 + \frac{1-\kappa^2}{\kappa^2}\sum_{\lambda \in m}\mathbb{E}(\hat{\beta}_\lambda - \beta_\lambda)^2 + \sum_{\lambda \in m}\mathbb{E}(\eta_{\lambda,\gamma}^2)\right\} + LD\sum_{\lambda \in \Gamma}F_\lambda. \tag{4.7}$$

In addition, there exists a constant $K_1$ depending on $p$, $\gamma$, $c$, $c'$, $\|f\|_1$ and on $\Phi$ such that

$$LD\sum_{\lambda \in \Gamma}F_\lambda \leq K_1(\log(n))^{c'+1}n^{c-\frac{\kappa^2\gamma}{q}-1}. \tag{4.8}$$

Since $\gamma > c$, one takes $\kappa < 1$ and $q > 1$ such that $c < \frac{\kappa^2\gamma}{q}$ and as required by Theorem 1, the last term satisfies

$$LD\sum_{\lambda \in \Gamma}F_\lambda \leq \frac{K_2}{n},$$

where $K_2$ is a constant. Before evaluating the first term, let us state the following lemma.

**Lemma 4.** *We set*

$$S_\varphi = \max\{\sup_{x \in supp(\phi)}|\phi(x)|, \sup_{x \in supp(\psi)}|\psi(x)|\}$$

*and*

$$I_\varphi = \min\{\inf_{x \in supp(\phi)}|\phi(x)|, \inf_{x \in supp(\psi)}|\psi(x)|\}.$$

*Using (2.4), we define $\Theta_\varphi = \frac{S_\varphi^2}{I_\varphi^2}$. For all $\lambda \in \Lambda$, we have the following result.*

- If $F_\lambda \leq \Theta_\varphi \frac{\log(n)}{n}$, then $\beta_\lambda^2 \leq \Theta_\varphi^2 \sigma_\lambda^2 \frac{\log(n)}{n}$.

- If $F_\lambda > \Theta_\varphi \frac{\log(n)}{n}$, then $\|\varphi_\lambda\|_\infty \frac{\log(n)}{n} \leq \sigma_\lambda \sqrt{\frac{\log(n)}{n}}$.

**Proof.** We note $\lambda = (j, k)$ and assume that $j \geq 0$ (arguments are similar for $j = -1$). If $F_\lambda \leq \Theta_\varphi \frac{\log(n)}{n}$, we have

$$|\beta_\lambda| \leq S_\varphi 2^{\frac{j}{2}} F_\lambda \leq S_\varphi 2^{\frac{j}{2}} \sqrt{F_\lambda} \sqrt{\Theta_\varphi} \sqrt{\frac{\log(n)}{n}} \leq S_\varphi I_\varphi^{-1} \sqrt{\Theta_\varphi} \sigma_\lambda \sqrt{\frac{\log(n)}{n}} \leq \Theta_\varphi \sigma_\lambda \sqrt{\frac{\log(n)}{n}},$$

since $\sigma_\lambda^2 \geq I_\varphi^2 2^j F_\lambda$. For the second point, observe that

$$\sigma_\lambda \sqrt{\frac{\log(n)}{n}} \geq 2^{\frac{j}{2}} I_\varphi \sqrt{\Theta_\varphi} \frac{\log(n)}{n} \text{ and } \|\psi_\lambda\|_\infty \frac{\log(n)}{n} \leq 2^{\frac{j}{2}} S_\varphi \frac{\log(n)}{n}.$$

∎

Now, for any $\delta > 0$,

$$\mathbb{E}(\eta_{\lambda,\gamma}^2) \leq (1+\delta) 2\gamma \log n \mathbb{E}(\tilde{V}_{\lambda,n}) + (1+\delta^{-1}) \left(\frac{\gamma \log n}{3n}\right)^2 \|\varphi_\lambda\|_\infty^2.$$

Moreover,

$$\mathbb{E}(\tilde{V}_{\lambda,n}) \leq (1+\delta) V_{\lambda,n} + (1+\delta^{-1}) 3\gamma \log n \frac{\|\varphi_\lambda\|_\infty^2}{n^2}.$$

So,

$$\mathbb{E}(\eta_{\lambda,\gamma}^2) \leq (1+\delta)^2 2\gamma \log n V_{\lambda,n} + \Delta(\delta) \left(\frac{\gamma \log n}{n}\right)^2 \|\varphi_\lambda\|_\infty^2, \tag{4.9}$$

with $\Delta(\delta)$ a constant depending only on $\delta$. Now, we apply (4.7) with

$$m = \left\{\lambda \in \Gamma_n : \beta_\lambda^2 > \Theta_\varphi^2 \frac{\sigma_\lambda^2}{n} \log n\right\},$$

so using Lemma 4, we can claim that for any $\lambda \in m$, $F_\lambda > \Theta_\varphi \frac{\log(n)}{n}$. Finally, since $\Theta_\varphi \geq 1$,

$$
\begin{aligned}
\mathbb{E}\|\tilde{\beta} - \beta\|_{\ell_2}^2 &\leq K_3 \left(\sum_{\lambda \in \Gamma_n} \beta_\lambda^2 1_{\{\beta_\lambda^2 \leq \Theta_\varphi^2 \frac{\sigma_\lambda^2}{n} \log n\}} + \sum_{\lambda \notin \Gamma_n} \beta_\lambda^2\right) \\
&\quad + K_3 \sum_{\lambda \in \Gamma_n} \left[\frac{\log n}{n} \sigma_\lambda^2 + \left(\frac{\log n}{n}\right)^2 \|\varphi_\lambda\|_\infty^2\right] 1_{\left\{\beta_\lambda^2 > \Theta_\varphi^2 \frac{\sigma_\lambda^2}{n} \log n, F_\lambda > \Theta_\varphi \frac{\log(n)}{n}\right\}} + \frac{K_4}{n} \\
&\leq K_3 \left[\sum_{\lambda \in \Gamma_n} \left(\beta_\lambda^2 1_{\{\beta_\lambda^2 \leq \Theta_\varphi^2 V_{\lambda,n} \log n\}} + 2\log n V_{\lambda,n} 1_{\{\beta_\lambda^2 > \Theta_\varphi^2 V_{\lambda,n} \log n\}}\right) + \sum_{\lambda \notin \Gamma_n} \beta_\lambda^2\right] + \frac{K_4}{n} \\
&\leq 2K_3 \left[\sum_{\lambda \in \Gamma_n} \min(\beta_\lambda^2, \Theta_\varphi^2 V_{\lambda,n} \log n) + \sum_{\lambda \notin \Gamma_n} \beta_\lambda^2\right] + \frac{K_4}{n},
\end{aligned}
$$

where the constant $K_3$ depends on $\gamma$ and $c$ and $K_4$ depends on $\gamma$, $c$, $c'$, $\|f\|_1$ and on $\Phi$. Theorem 1 is proved by using properties of the biorthogonal wavelet basis.

### 4.4 Proof of Theorem 2

Let us assume that $f$ belongs to $B_{2,\Gamma}^s(R^{1-2s}) \cap W_s(R) \cap \mathbb{L}_1(R) \cap \mathbb{L}_2(R)$. Inequality (1.6) of Theorem 1 implies that, for all $n$,

$$\mathbb{E}(\|\tilde{f}_{n,\gamma} - f\|_2^2) \leq C_1 \left[ \sum_{\lambda \in \Gamma_n} \left( \beta_\lambda^2 1_{|\beta_\lambda| \leq \sigma_\lambda \sqrt{\frac{\log n}{n}}} + V_{\lambda,n} \log n 1_{|\beta_\lambda| > \sigma_\lambda \sqrt{\frac{\log n}{n}}} \right) + \sum_{\lambda \notin \Gamma_n} \beta_\lambda^2 \right] + \frac{C_2}{n}$$

where $C_1$ and $C_2$ are two constants. But we have

$$\sum_{\lambda \in \Gamma_n} V_{\lambda,n} \log n 1_{|\beta_\lambda| > \sigma_\lambda \sqrt{\frac{\log n}{n}}} = \sum_{\lambda \in \Gamma_n} \sigma_\lambda^2 \frac{\log n}{n} \sum_{k=0}^{+\infty} 1_{2^{-k-1}\beta_\lambda^2 \leq \sigma_\lambda^2 \frac{\log n}{n} < 2^{-k}\beta_\lambda^2}$$

$$\leq \sum_{k=0}^{+\infty} 2^{-k} \sum_{\lambda \in \Lambda} \beta_\lambda^2 1_{|\beta_\lambda| \leq 2^{\frac{k+1}{2}} \sigma_\lambda \sqrt{\frac{\log n}{n}}}$$

$$\leq \sum_{k=0}^{+\infty} 2^{-k} R^{2-4s} \left( 2^{\frac{k+1}{2}} \sqrt{\frac{\log n}{n}} \right)^{4s}$$

$$\leq R^{2-4s} \rho_{n,s}^2 \sum_{k=0}^{+\infty} 2^{-k+2s(k+1)}$$

and

$$\sum_{\lambda \notin \Gamma_n} \beta_\lambda^2 \leq R^{2-4s} \rho_{n,s}^2.$$

So,

$$\mathbb{E}(\|\tilde{f}_{n,\gamma} - f\|_2^2) \leq C(\gamma, c, \Phi, s) R^{2-4s} \rho_{n,s}^2 + \frac{C_2}{n},$$

where $C(\gamma, c, \Phi, s)$ depends on $\gamma$, $c$, $\Phi$ and $s$. Hence,

$$\mathbb{E}(\|\tilde{f}_{n,\gamma} - f\|_2^2) \leq C(\gamma, c, \Phi, s) R^{2-4s} \rho_{n,s}^2 (1 + o_n(1))$$

and $f$ belongs to $MS(\tilde{f}_\gamma, \rho_s)(R')$ for $R'$ large enough.

Conversely, let us suppose that $f$ belongs to $MS(\tilde{f}_\gamma, \rho_s)(R') \cap \mathbb{L}_1(R') \cap \mathbb{L}_2(R')$. Then, for any $n$,

$$\mathbb{E}(\|\tilde{f}_{n,\gamma} - f\|_2^2) \leq R'^2 \left( \frac{\log n}{n} \right)^{2s}.$$

Consequently, there exists $R$ depending on $R'$ and $\Phi$ such that for any $n$,

$$\sum_{\lambda \notin \Gamma_n} \beta_\lambda^2 \leq R^2 \left( \frac{\log n}{n} \right)^{2s}.$$

This implies that $f$ belongs to $B_{2,\Gamma}^s(R)$.

Now, we want to prove that $f \in W_s(R)$ if $R$ is large enough. We have

$$\sum_{\lambda \in \Lambda} \beta_\lambda^2 1_{|\beta_\lambda| \leq \sigma_\lambda \sqrt{\frac{\gamma \log n}{2n}}} \leq \sum_{\lambda \notin \Gamma_n} \beta_\lambda^2 + \sum_{\lambda \in \Gamma_n} \beta_\lambda^2 1_{|\beta_\lambda| \leq \sigma_\lambda \sqrt{\frac{\gamma \log n}{2n}}}.$$

But $\tilde{\beta}_\lambda = \hat{\beta}_\lambda 1_{|\hat{\beta}_\lambda| \geq \eta_{\lambda,\gamma}}$, so,

$$|\beta_\lambda| 1_{|\beta_\lambda| \leq \frac{\eta_{\lambda,\gamma}}{2}} \leq |\beta_\lambda - \tilde{\beta}_\lambda|.$$

So, for any $n$,

$$\begin{aligned}
\sum_{\lambda \in \Lambda} \beta_\lambda^2 1_{|\beta_\lambda| \leq \sigma_\lambda \sqrt{\frac{\gamma \log n}{2n}}} &\leq \sum_{\lambda \notin \Gamma_n} \beta_\lambda^2 + \mathbb{E}\left\{ \sum_{\lambda \in \Gamma_n} \beta_\lambda^2 1_{|\beta_\lambda| \leq \sigma_\lambda \sqrt{\frac{\gamma \log n}{2n}}} [1_{|\beta_\lambda| \leq \frac{\eta_{\lambda,\gamma}}{2}} + 1_{|\beta_\lambda| > \frac{\eta_{\lambda,\gamma}}{2}}] \right\} \\
&\leq \sum_{\lambda \notin \Gamma_n} \beta_\lambda^2 + \sum_{\lambda \in \Gamma_n} \mathbb{E}[(\tilde{\beta}_\lambda - \beta_\lambda)^2] + \sum_{\lambda \in \Gamma_n} \beta_\lambda^2 1_{|\beta_\lambda| \leq \sigma_\lambda \sqrt{\frac{\gamma \log n}{2n}}} \mathbb{E}(1_{|\hat{\beta}_\lambda| > \frac{\eta_{\lambda,\gamma}}{2}}) \\
&\leq \sum_{\lambda \notin \Gamma_n} \beta_\lambda^2 + \sum_{\lambda \in \Gamma_n} \mathbb{E}[(\tilde{\beta}_\lambda - \beta_\lambda)^2] + \sum_{\lambda \in \Gamma_n} \beta_\lambda^2 \mathbb{P}\left( \sigma_\lambda \sqrt{\frac{\gamma \log n}{2n}} > \frac{\eta_{\lambda,\gamma}}{2} \right) \\
&\leq \mathbb{E}(\|\tilde{\beta} - \beta\|_{\ell_2}^2) + \sum_{\lambda \in \Gamma_n} \beta_\lambda^2 \mathbb{P}\left( \sigma_\lambda \sqrt{\frac{\gamma \log n}{2n}} > \frac{\eta_{\lambda,\gamma}}{2} \right).
\end{aligned}$$

Using Lemma 1,

$$\mathbb{P}\left( \sigma_\lambda \sqrt{\frac{2\gamma \log n}{n}} > \eta_{\lambda,\gamma} \right) \leq \mathbb{P}(\tilde{V}_{\lambda,n} \leq V_{\lambda,n}) \leq n^{-\gamma}$$

and

$$\sum_{\lambda \in \Lambda} \beta_\lambda^2 1_{|\beta_\lambda| \leq \sigma_\lambda \sqrt{\frac{\gamma \log n}{2n}}} \leq c_1(\Phi)^{-1}(R')^2 \left( \sqrt{\frac{\log n}{n}} \right)^{4s} + \|\beta\|_{\ell_2}^2 n^{-\gamma}.$$

Since this is true for every $n$, we have for any $t \leq 1$,

$$\sum_{\lambda \in \Lambda} \beta_\lambda^2 1_{|\beta_\lambda| \leq \sigma_\lambda t} \leq R^{2-4s} \left( \sqrt{\frac{2}{\gamma}} t \right)^{4s},$$

where $R$ is a constant large enough depending on $R'$ and $\Phi$. Note that

$$\sup_{t \geq 1} t^{-4s} \sum_{\lambda \in \Lambda} \beta_\lambda^2 1_{|\beta_\lambda| \leq \sigma_\lambda t} \leq \|\beta\|_{\ell_2}^2.$$

We conclude that

$$f \in B_{2,\Gamma}^s(R) \cap W_s(R)$$

for $R$ large enough.

### 4.5 Proof of Proposition 2

Since $\beta < \frac{1}{2}$, $f_\beta \in \mathbb{L}_1 \cap \mathbb{L}_2$. If the Haar basis is considered, the wavelet coefficients $\beta_{j,k}$ of $f_\beta$ can be calculated and we obtain for any $j \geq 0$, for any $k \notin \{0, \ldots, 2^j - 1\}$, $\beta_{j,k} = 0$ and for any $j \geq 0$, for any $k \in \{0, \ldots, 2^j - 1\}$,

$$\beta_{j,k} = (1-\beta)^{-1} 2^{-j(\frac{1}{2}-\beta)} \left( 2\left(k + \frac{1}{2}\right)^{1-\beta} - k^{1-\beta} - (k+1)^{1-\beta} \right)$$

and there exists a constant $0 < c_{1,\beta} < \infty$ only depending on $\beta$ such that

$$\lim_{k \to \infty} 2^{j(\frac{1}{2}-\beta)} k^{1+\beta} \beta_{j,k} = c_{1,\beta}.$$

Moreover the $\beta_{j,k}$'s are strictly positive. Consequently they can be upper and lower bounded, up to a constant, by $2^{-j(\frac{1}{2}-\beta)} k^{-(1+\beta)}$. Similarly, for any $j \geq 0$, for any $k \in \{0, \ldots, 2^j - 1\}$,

$$\sigma_{j,k}^2 = (1-\beta)^{-1} 2^{j\beta} \left((k+1)^{1-\beta} - k^{1-\beta}\right)$$

and there exists a constant $0 < c_{2,\beta} < \infty$ only depending on $\beta$ such that

$$\lim_{k \to \infty} 2^{-j\beta} k^\beta \sigma_{j,k}^2 = c_{2,\beta}.$$

There exist two constants $\kappa(\beta)$ and $\kappa'(\beta)$ only depending on $\beta$ such that for any $0 < t < 1$, if $\beta_{j,k} \neq 0$

$$|\beta_{j,k}| \leq t\sigma_{j,k} \Rightarrow k \geq \kappa(\beta) t^{-\frac{2}{\beta+2}} 2^{j\left(\frac{\beta-1}{\beta+2}\right)}$$

and

$$\kappa(\beta) t^{-\frac{2}{\beta+2}} 2^{j\left(\frac{\beta-1}{\beta+2}\right)} \geq 2^j \iff 2^j \leq \kappa'(\beta) t^{-\frac{2}{3}}.$$

So, if $2^j \leq \kappa'(\beta) t^{-\frac{2}{3}}$, since $\beta_{jk} = 0$ for $k \geq 2^j$,

$$\sum_{k \in \mathbb{Z}} \beta_{j,k}^2 \mathbf{1}_{\beta_{j,k} \leq t\sigma_{j,k}} = 0.$$

We obtain

$$\sum_{\lambda \in \Lambda} \beta_\lambda^2 \mathbf{1}_{|\beta_\lambda| \leq t\sigma_\lambda} \leq C(\beta) \sum_{j=-1}^{+\infty} 2^{-j(1-2\beta)} \mathbf{1}_{2^j > \kappa'(\beta) t^{-\frac{2}{3}}} \sum_{k=1}^{2^j-1} k^{-2-2\beta} \leq C'(\beta) t^{\frac{2-4\beta}{3}},$$

where $C(\beta)$ and $C'(\beta)$ denote two constants only depending on $\beta$. So, for any $0 < s < \frac{1}{6}$, if we take $\beta \leq \frac{1}{2}(1-6s)$, then, for any $0 < t < 1$, $t^{\frac{2-4\beta}{3}} \leq t^{4s}$. Finally, there exists $c \geq 1$, such that for any $n$,

$$\sum_{\lambda \notin \Gamma_n} \beta_\lambda^2 \leq R^2 \rho_{n,\alpha}^2,$$

where $R > 0$. And in this case,

$$f_\beta \notin \mathbb{L}_\infty, \quad f_\beta \in B_{2,\infty}^s \cap W_s := MS(\tilde{f}_\gamma^H, \rho_s).$$

### 4.6 Proof of Theorem 3

Using the maxiset results of Section 3.1, since

$$MS(\tilde{f}_\gamma, \rho_{\frac{\alpha}{1+2\alpha}}) :=: \mathcal{B}_{2,\infty}^{\frac{\alpha}{c(1+2\alpha)}} \cap W_{\frac{\alpha}{1+2\alpha}},$$

it is enough to show that

$$\mathcal{B}_{p,q}^\alpha(R) \cap \mathcal{L}_{1,2,\infty}(R') \subset \mathcal{B}_{2,\infty}^{\frac{\alpha}{c(1+2\alpha)}}(R'') \cap W_{\frac{\alpha}{1+2\alpha}}(R'')$$

for $R'' > 0$ (see (3.1)). Let $f \in \mathcal{B}^\alpha_{p,q}(R) \cap \mathcal{L}_{1,2,\infty}(R')$. We first prove that $f \in W_{\frac{\alpha}{1+2\alpha}}(R'')$ for $R''$ large enough. Since for any $\lambda = (j,k)$,

$$\sigma_\lambda^2 \leq \min\left(\max(2^j;1)\|\varphi\|_\infty^2 F_{j,k}\, ; \, \|f\|_\infty \|\varphi\|_2^2\right),$$

where $\varphi \in \{\phi, \psi\}$ according to the value of $j$, we have for any $t > 0$ and any $\tilde{J}$

$$\begin{aligned}
\sum_\lambda \beta_\lambda^2 1_{|\beta_\lambda| \leq \sigma_\lambda t} &\leq \sum_{j<\tilde{J}} \sum_k \sigma_{j,k}^2 t^2 + \sum_{j \geq \tilde{J}} \sum_k \beta_{j,k}^2 \left(\frac{\sigma_{j,k} t}{|\beta_{j,k}|}\right)^{2-p} \\
&\leq \max(\|\phi\|_\infty^2; \|\psi\|_\infty^2) t^2 \sum_{j<\tilde{J}} \max(2^j;1) \sum_k F_{j,k} + \sum_{j \geq \tilde{J}} \sum_k \beta_{j,k}^2 \left(\frac{t\sqrt{\|f\|_\infty}\|\psi\|_2^2}{|\beta_{j,k}|}\right)^{2-p} \\
&\leq C(\Phi, R') \left(2^{\tilde{J}} t^2 + t^{2-p} \sum_{j \geq \tilde{J}} \sum_k |\beta_{j,k}|^p\right),
\end{aligned}$$

where $C(\Phi, R')$ is a constant only depending on $\Phi$ and on $R'$. Indeed, we have used that

$$\sum_k F_{j,k} \leq m_\varphi \|f\|_1, \tag{4.10}$$

by similar arguments to (4.2)). Now, since $f$ belongs to $\mathcal{B}^\alpha_{p,\infty}(R)$ (that contains $\mathcal{B}^\alpha_{p,q}(R)$, see Section 2.2), with $\alpha + \frac{1}{2} - \frac{1}{p} > 0$,

$$\sum_\lambda \beta_\lambda^2 1_{|\beta_\lambda| \leq \sigma_\lambda t} \leq C_1(\Phi, \alpha, p, R') \left(2^{\tilde{J}} t^2 + t^{2-p} R^p 2^{-\tilde{J} p(\alpha + \frac{1}{2} - \frac{1}{p})}\right),$$

where $C_1(\Phi, \alpha, p, R')$ depends on $\Phi$, $\alpha$, $p$ and $R'$. With $\tilde{J}$ such that

$$2^{\tilde{J}} \leq R^{\frac{2}{1+2\alpha}} t^{-\frac{2}{1+2\alpha}} < 2^{\tilde{J}+1},$$

$$\sum_\lambda \beta_\lambda^2 1_{|\beta_\lambda| \leq \sigma_\lambda t} \leq C_2(\Phi, \alpha, p, R') R^{\frac{2}{1+2\alpha}} t^{\frac{4\alpha}{1+2\alpha}}$$

where $C_2(\Phi, \alpha, p, R')$ depends on $\Phi$, $\alpha$, $p$ and $R'$. So, $f$ belongs to $W_{\frac{\alpha}{1+2\alpha}}(R'')$ for $R''$ large enough. Furthermore, using (2.5), if $p \leq 2$ and

$$\alpha\left(1 - \frac{1}{c(1+2\alpha)}\right) \geq \frac{1}{p} - \frac{1}{2}$$

$$\mathcal{B}^\alpha_{p,\infty}(R) \subset \mathcal{B}^{\frac{\alpha}{c(1+2\alpha)}}_{2,\infty}(R).$$

Finally, for $R''$ large enough,

$$\mathcal{B}^\alpha_{p,q}(R) \cap \mathcal{L}_{1,2,\infty}(R') \subset \mathcal{B}^\alpha_{p,\infty}(R) \cap \mathcal{L}_{1,2,\infty}(R') \subset \mathcal{B}^{\frac{\alpha}{c(1+2\alpha)}}_{2,\infty}(R'') \cap W_{\frac{\alpha}{1+2\alpha}}(R'').$$

### 4.7  Proof of Theorem 4

In this subsection since $\alpha > 0$ and $p > 2$, we set

$$s = \frac{\alpha}{2\alpha + 2 - \frac{2}{p}}.$$

Using the maxiset results of Section 3.1, since

$$MS(\tilde{f}_\gamma, \rho_s) :=: \mathcal{B}_{2,\infty}^{c^{-1}s} \cap W_s,$$

it is enough to show that

$$\mathcal{B}_{p,q}^{\alpha}(R) \cap \mathbb{L}_1(R') \cap \mathbb{L}_2(R') \subset \mathcal{B}_{2,\infty}^{c^{-1}s}(R'') \cap W_s(R'')$$

for $R'' > 0$ (see (3.1)). By using (2.5), since $c \geq 1$, we have

$$\mathcal{B}_{p,q}^{\alpha}(R) \subset \mathcal{B}_{p,\infty}^{\alpha}(R) \subset \mathcal{B}_{2,\infty}^{c^{-1}s}(R).$$

Let $f \in \mathcal{B}_{p,q}^{\alpha}(R) \cap \mathbb{L}_1(R') \cap \mathbb{L}_2(R')$. We prove that $f \in W_s(R'')$ for $R''$ large enough. Using computations of Section 4.6, we have for any $t > 0$ and any $\tilde{J} \geq 0$

$$\sum_\lambda \beta_\lambda^2 \mathbf{1}_{|\beta_\lambda| \leq \sigma_\lambda t} \leq C(\Phi, R') \left( 2^{\tilde{J}} t^2 + \sum_{j \geq \tilde{J}} \sum_k \beta_{j,k}^2 \right),$$

where $C(\Phi, R')$ is a constant only depending on $\Phi$ and on $R'$. Now, let us bound for all $j \geq \tilde{J}$

$$\sum_k \beta_{j,k}^2 = \sum_k |\beta_{j,k}|^{\frac{p}{p-1}} |\beta_{j,k}|^{2-\frac{p}{p-1}}.$$

Let us apply the Hölder inequality. Since $p > 2$, we have $2 - \frac{p}{p-1} > 0$ and

$$\sum_k \beta_{j,k}^2 \leq \left( \sum_k |\beta_{j,k}|^p \right)^{\frac{1}{p-1}} \left( \sum_k |\beta_{j,k}| \right)^{2-\frac{p}{p-1}}.$$

Since $f \in \mathcal{B}_{p,\infty}^{\alpha}(R)$,

$$\left( \sum_k |\beta_{j,k}|^p \right)^{\frac{1}{p-1}} \leq R^{\frac{p}{p-1}} 2^{-\frac{jp}{p-1}(\alpha + \frac{1}{2} - \frac{1}{p})}.$$

Since $f \in \mathbb{L}_1(R')$,

$$\sum_k |\beta_{j,k}| = \sum_k \left| 2^{\frac{j}{2}} \int f(x) \psi(2^j x - k) dx \right|$$

$$\leq 2^{\frac{j}{2}} \|\psi\|_\infty \sum_k F_{jk}$$

$$\leq 2^{\frac{j}{2}} \|\psi\|_\infty m_\varphi \|f\|_1$$

by using (4.10). Hence

$$\sum_k \beta_{j,k}^2 \leq R^{\frac{p}{p-1}} (\|\psi\|_\infty m_\varphi R')^{2-\frac{p}{p-1}} 2^{-j\alpha \frac{p}{p-1}}.$$

Finally,
$$\sum_\lambda \beta_\lambda^2 1_{|\beta_\lambda| \leq \sigma_\lambda t} \leq C_1(\Phi, R') \left(2^{\tilde{J}} t^2 + R^{\frac{p}{p-1}} 2^{-\tilde{J}\alpha \frac{p}{p-1}}\right)$$

where $C_1(\Phi, R')$ is a constant only depending on $\Phi$ and on $R'$. With $\tilde{J}$ such that

$$2^{\tilde{J}} \leq R^{\frac{p}{\alpha p+p-1}} t^{-\frac{2(p-1)}{\alpha p+p-1}} < 2^{\tilde{J}+1},$$

$$\sum_\lambda \beta_\lambda^2 1_{|\beta_\lambda| \leq \sigma_\lambda t} \leq C_2(\Phi, \alpha, p, R') R^{\frac{1}{\alpha+1-\frac{1}{p}}} t^{\frac{2\alpha}{\alpha+1-\frac{1}{p}}}$$

where $C_2(\Phi, \alpha, p, R')$ depends on $\Phi$, $\alpha$, $p$ and $R'$. So, $f$ belongs to $W_s(R'')$ for $R''$ large enough. Finally, for $R''$ large enough,

$$\mathcal{B}_{p,q}^\alpha(R) \cap \mathbb{L}_1(R') \cap \mathbb{L}_2(R') \subset \mathcal{B}_{2,\infty}^{c^{-1}s}(R'') \cap W_s(R'').$$

### 4.8 Proof of Theorem 5

To establish the lower bound stated in Theorem 5, we first consider $p \geq 2$ and $0 < \alpha < r+1$. As usual, the lower bound of the risk

$$\mathcal{R}_n(\alpha, p) = \inf_{\hat{f}} \sup_{f \in \mathcal{B}_{p,\infty}^\alpha(R) \cap \mathbb{L}_1(R_1) \cap \mathbb{L}_2(R_2) \cap \mathbb{L}_\infty(R_\infty)} \mathbb{E}\left[\|f - \hat{f}\|_2^2\right],$$

where $R$, $R_1$, $R_2$ and $R_\infty$ are positive real numbers, can be obtained by using an adequate version of Fano's lemma based on the Kullback-Leibler divergence. We first give classical lemmas that introduce constants useful in the sequel. The first result recalls the Kullback-Leibler divergence for Poisson processes (see [10]).

**Lemma 5.** *Let $N$ and $N'$ be two Poisson processes on $\mathbb{R}$ whose intensities with respect to the Lebesgue measure are respectively $s$ and $s'$. We denote $\mathbb{P}$ (respectively $\mathbb{Q}$) the probability measures associated with $s$ (respectively with $s'$). Then, the Kullback-Leibler divergence between $\mathbb{P}$ and $\mathbb{Q}$ is*

$$K(\mathbb{P}, \mathbb{Q}) = \int_\mathbb{X} s(x) \phi\left(\log\left(\frac{s'(x)}{s(x)}\right)\right) dx$$

*where $\phi(u) = \exp(u) - u - 1$.*

Now, let us give the following version of Fano's lemma, derived from [6].

**Lemma 6.** *Let $(\mathbb{P}_i)_{i \in \{0,...,n\}}$ be a finite family of probability measures defined on the same measurable space $\Omega$. One sets*

$$\bar{K}_n = \frac{1}{n} \sum_{i=1}^n K(\mathbb{P}_i, \mathbb{P}_0).$$

*Then, there exists an absolute constant $B$ ($B = 0.71$ works) such that if $\hat{\theta}$ is a random variable on $\Omega$ with values in $\{0,...,n\}$, one has*

$$\inf_{0 \leq i \leq n} \mathbb{P}_i(\hat{\theta} = i) \leq \max\left(B, \frac{\bar{K}_n}{\log(n+1)}\right).$$

Finally, we recall a combinatorial lemma due to Gallager (see Lemma 8 in [30]).

**Lemma 7.** *Let $\Gamma$ be a finite set with cardinal $Q$. Let $D \leq Q$. There exist absolute constants $\theta$ and $\sigma$ such that there exists $\mathcal{M}_D \subset \mathcal{P}(\Gamma)$, satisfying $\log|\mathcal{M}_D| \geq \sigma D$ if $D = Q$ and $\log|\mathcal{M}_D| \geq \sigma D \log(Q/D)$ if $D < Q$ and such that for all distinct sets $m$ and $m'$ belonging to $\mathcal{M}_D$ we have $|m \triangle m'| \geq \theta D$.*

Now, we are ready to provide a lower bound for $\mathcal{R}_n(\alpha, p)$. For this purpose, for a given $n$ large enough, we set $j$ the largest integer such that

$$2^j \leq \left(\frac{R}{2B\sigma c_2(\Phi)^{-1}c_{\tilde{\psi}}}\right)^{\frac{1}{\alpha+1-\frac{1}{p}}} \left(\frac{R_1}{2B\sigma c_2(\Phi)^{-1}c_{\tilde{\psi}}^2}\right)^{\frac{-1}{p\alpha+p-1}} n^{\frac{1-\frac{1}{p}}{\alpha+1-\frac{1}{p}}}.$$

The constant $c_2(\Phi)$ was defined in Section 2.2 and $c_{\tilde{\psi}}$ is a constant depending only on $\tilde{\psi}$ such that

$$\|\sum_{k \in \mathbb{Z}} \tilde{\psi}_{0,k}\|_\infty \leq c_{\tilde{\psi}}.$$

We set for any $\ell$,

$$g_\ell(x) = \frac{\int_0^{x-\ell+1} \exp\left(-\frac{1}{u(1-u)}\right) du}{\int_0^1 \exp\left(-\frac{1}{u(1-u)}\right) du} 1_{[\ell-1,\ell]}(x) + 1_{]\ell,\ell+1]}(x).$$

Note that $\delta = \|g_\ell\|_1$ does not depend on $\ell$. We also introduce the integer $D$ such that $D2^{-j}$ is the largest integer satisfying

$$D2^{-j} \leq \frac{R_1 n 2^{-j}}{2B\sigma c_2(\Phi)^{-1}c_{\tilde{\psi}}^2} - 2\delta. \tag{4.11}$$

In particular, $D2^{-j}$ goes to $\infty$ when $n$ goes to $\infty$. Using Lemma 7 with $\Gamma = \{0, 1, \ldots, D-1\}$ and $Q = D$, we extract $\mathcal{M}_D$ for which both properties stated in Lemma 7 are satisfied and we set

$$\mathcal{C}_{j,D} = \left\{f_m = \tilde{f}_{j,D} + a_j \sum_{k \in m} \tilde{\psi}_{j,k} : \quad m \in \mathcal{M}_D\right\},$$

with

$$a_j = \frac{B\sigma c_2(\Phi)^{-1} c_{\tilde{\psi}} 2^{\frac{j}{2}}}{n}.$$

The function $\tilde{f}_{j,D}$ is defined by

$$\tilde{f}_{j,D}(x) = \rho 1_{[0,D2^{-j}]}(x) + \rho g_{-1}(x) + \rho g_{-D2^{-j}-1}(-x)$$

where

$$\rho = \frac{R_1 2^j D^{-1}}{1 + 2\delta 2^j D^{-1}}.$$

Let $f_m \in \mathcal{C}_{j,D}$. Observe that the support of $\sum_{k \in m} \tilde{\psi}_{j,k}$ is included in $[-1, D2^{-j}+1]$ for $n$ large enough. In this case, since $\rho \geq 2a_j 2^{\frac{j}{2}} c_{\tilde{\psi}}$ (see (4.11)), we have for $x$ in the support of $\sum_{k \in m} \tilde{\psi}_{j,k}$

$$f_m(x) \geq \frac{\rho}{2}. \tag{4.12}$$

In addition for any $x$, $f_m(x) \geq 0$. Now, we verify that $f_m$ belongs to $\mathcal{B}^\alpha_{p,\infty}(R) \cap \mathbb{L}_1(R_1) \cap \mathbb{L}_2(R_2) \cap \mathbb{L}_\infty(R_\infty)$. We have:

$$
\begin{aligned}
\|f_m\|_{\alpha,p,\infty} &\leq \|\tilde{f}_{j,D}\|_{\alpha,p,\infty} + \|a_j \sum_{k \in m} \tilde{\psi}_{j,k}\|_{\alpha,p,\infty} \\
&\leq \|\tilde{f}_{j,D}\|_{\alpha,p,\infty} + D^{\frac{1}{p}} a_j 2^{j(\alpha+\frac{1}{2}-\frac{1}{p})} \\
&\leq \|\tilde{f}_{j,D}\|_{\alpha,p,\infty} + \left(\frac{R_1 n 2^{-j}}{2B\sigma c_2(\Phi)^{-1} c_{\tilde{\psi}}^2}\right)^{\frac{1}{p}} \frac{B\sigma c_2(\Phi)^{-1} c_{\tilde{\psi}} 2^{\frac{j}{2}}}{n} 2^{j(\alpha+\frac{1}{2})} \\
&= \|\tilde{f}_{j,D}\|_{\alpha,p,\infty} + 2^{j(\alpha+1-\frac{1}{p})} \left(\frac{R_1}{2B\sigma c_2(\Phi)^{-1} c_{\tilde{\psi}}^2}\right)^{\frac{1}{p}} B\sigma c_2(\Phi)^{-1} c_{\tilde{\psi}} n^{\frac{1}{p}-1} \\
&\leq \|\tilde{f}_{j,D}\|_{\alpha,p,\infty} + \frac{R}{2}.
\end{aligned}
$$

Finally, $\tilde{f}_{j,D}$ has an infinite number of continuous derivatives bounded (up to constants) by $\rho$ and $\|\tilde{f}_{j,D}\|_{\alpha,p,\infty}$ is bounded (up to a constant) by $\rho(D2^{-j})^{1/p}$ that goes to 0 when $n$ goes to $\infty$. So, for $n$ large enough,

$$\|f_m\|_{\alpha,p,\infty} \leq R.$$

Now, it remains to verify that $f_m \in \mathbb{L}_1(R_1) \cap \mathbb{L}_2(R_2) \cap \mathbb{L}_\infty(R_\infty)$. We have

$$\|f_m\|_\infty \leq \rho + c_{\tilde{\psi}} 2^{\frac{j}{2}} a_j \leq R_1 2^j D^{-1} + \frac{B\sigma c_2(\Phi)^{-1} c_{\tilde{\psi}}^2 2^j}{n} \leq R_\infty$$

for $n$ large enough. Using again (4.11),

$$\|f_m\|_2^2 \leq 2\|\tilde{f}_{j,D}\|_2^2 + 2\|a_j \sum_{k \in m} \tilde{\psi}_{j,k}\|_2^2 \leq 2\rho^2(D2^{-j} + 2\delta) + 2c_2(\Phi)Da_j^2 \leq 2\rho R_1 + \frac{R_1 B\sigma 2^j}{n} \leq R_2^2$$

for $n$ large enough. Since $f_m \geq 0$,

$$\|f_m\|_1 = \int_{-\infty}^{+\infty} \left(\tilde{f}_{j,D}(x) + a_j \sum_{k \in m} \tilde{\psi}_{j,k}(x)\right) dx = \rho D 2^{-j} + 2\delta\rho = R_1.$$

Finally, we have:
$$\mathcal{R}_n(\alpha,p) \geq \inf_{\hat{f}} \sup_{f \in \mathcal{C}_{j,D}} \mathbb{E}\left[\|f - \hat{f}\|_2^2\right].$$

If $\hat{f}$ is an estimator, we can define $\hat{f}' = \arg\min_{t \in \mathcal{C}_{j,D}} \|t - \hat{f}\|_2$. Then, for $f \in \mathcal{C}_{j,D}$,

$$\|\hat{f}' - f\|_2 \leq \|\hat{f}' - \hat{f}\|_2 + \|\hat{f} - f\|_2 \leq 2\|\hat{f} - f\|_2$$

and

$$\mathcal{R}_n(\alpha,p) \geq \frac{1}{4} \inf_{\hat{f} \in \mathcal{C}_{j,D}} \sup_{f \in \mathcal{C}_{j,D}} \mathbb{E}\left[\|f - \hat{f}\|_2^2\right].$$

Moreover if $m$ and $m'$ belong to $\mathcal{M}_D$ with $m \neq m'$,

$$\|f_m - f_{m'}\|_2^2 \geq c_1(\Phi) a_j^2 |m \Delta m'| \geq c_1(\Phi) \theta D a_j^2$$

where $c_1(\Phi)$ has been defined in Section 2.2. Hence

$$\mathcal{R}_n(\alpha, p) \geq \frac{c_1(\Phi)}{4} \theta D a_j^2 \inf_{\hat{f} \in \mathcal{C}_D} \sup_{f \in \mathcal{C}_D} \mathbb{P}_f(\hat{f} \neq f).$$

To apply Lemma 6, we need to compute $\bar{K}_n$. For any distinct sets $m$ and $m'$ belonging to $\mathcal{M}_D$, since for any $x > -1$, $\log(1 + x) \geq x/(1+x)$ and by using (4.12), we have

$$\begin{aligned}
K(\mathbb{P}_{f_{m'}}, \mathbb{P}_{f_m}) &= \int f_{m'} \phi(\log \frac{f_m}{f_{m'}}) n dx \\
&= \int [f_m - f_{m'} - f_{m'} \log(1 + \frac{f_m - f_{m'}}{f_{m'}})] n dx \\
&\leq \int \frac{(f_m - f_{m'})^2}{f_m}(x) n dx \\
&\leq \frac{2}{\rho} n \|f_{m'} - f_m\|_2^2 \\
&\leq \frac{2 n a_j^2 D c_2(\Phi)}{\rho}
\end{aligned} \qquad (4.13)$$

and $\bar{K}_n \leq \frac{2na_j^2 D c_2(\Phi)}{\rho}$. By applying Lemma 6, since

$$\frac{2c_2(\Phi) n D a_j^2}{\rho \sigma D} \leq B,$$

we have

$$\begin{aligned}
\mathcal{R}_n(\alpha, p) &\geq \frac{c_1(\Phi)}{4} \theta (1 - B) D a_j^2 \\
&\geq \frac{c_1(\Phi)}{4} \theta (1 - B) \frac{R_1 n}{2 B \sigma c_2(\Phi)^{-1} c_{\tilde{\psi}}^2} \frac{(B \sigma c_2(\Phi)^{-1} c_{\tilde{\psi}})^2 2^j}{n^2} (1 + o_n(1)) \\
&\geq C R^{\frac{1}{\alpha + 1 - \frac{1}{p}}} n^{-\frac{\alpha}{\alpha + 1 - \frac{1}{p}}} (1 + o_n(1)),
\end{aligned}$$

where $C$ is a constant that depends on $\alpha$, $p$, $c_2(\Phi)$, $c_{\tilde{\psi}}$, $\theta$, $B$, $\sigma$ and $R_1$, which is the stated result.

For the case $p \leq 2$, by using computations similar to those of Theorem 2 of [18], it is easy to prove that the minimax risk associated to the set of functions supported by $[0,1]$ and belonging to $\mathcal{B}_{p,q}^\alpha(R)$ for $0 < \alpha < r + 1$ is larger than $n^{-\frac{2\alpha}{1+2\alpha}}$ up to a constant.

Finally, the adaptive properties of $\tilde{f}_\gamma$ are proved by combining Theorems 3 and 4 and the previous lower bound.

### 4.9 Proof of Theorem 6

Let us consider the Haar basis. For $j \geq 0$ and $D \in \{0, 1, \ldots, 2^j\}$, we set

$$\mathcal{C}_{j,D} = \{f_m = \rho 1_{[0,1]} + a_{j,D} \sum_{k \in m} \tilde{\varphi}_{j,k} : \quad |m| = D, m \subset \mathcal{N}_j\},$$

where
$$\mathcal{N}_j = \{k : \quad \tilde{\varphi}_{j,k} \text{ has support in } [0,1]\}.$$

The parameters $j, D, \rho, a_{j,D}$ is chosen later to fulfill some requirements. Note that
$$N_j = \operatorname{card}(\mathcal{N}_j) = 2^j.$$

We know that there exists a subset of $\mathcal{C}_{j,D}$, denoted $\mathcal{M}_{j,D}$, and some universal constants, denoted $\theta$ and $\sigma$, such that for all $m, m' \in \mathcal{M}_{j,D}$,
$$\operatorname{card}(m \Delta m') \geq \theta D, \quad \log(\operatorname{card}(\mathcal{M}_{j,D})) \geq \sigma D \log\left(\frac{2^j}{D}\right)$$

(see Lemma 7). Now, let us describe all the requirements necessary to obtain the lower bound of the risk.

- To ensure $f_m \geq 0$ and the equivalence between the Kullback distance and the $\mathbb{L}_2$-norm (see below), the $f_m$'s have to be larger than $\rho/2$. Since the $\tilde{\varphi}_{j,k}$'s have disjoint support, this means that
$$\rho \geq 2^{1+j/2}|a_{j,D}|. \tag{4.14}$$

- We need the $f_m$'s to be in $\mathbb{L}_1(R'') \cap \mathbb{L}_\infty(R'')$. Since $\|f\|_1 = \rho$ and $\|f\|_\infty = \rho + 2^{j/2}|a_{j,D}|$, we need
$$\rho + 2^{j/2}|a_{j,D}| \leq R''. \tag{4.15}$$

- The $f_m$'s have to belong to $B_{2,\Gamma}^s(R')$ i.e.
$$\rho + 2^{js}\sqrt{D}|a_{j,D}| \leq R'. \tag{4.16}$$

- The $f_m$'s have to belong to $W_s(R)$. We have $\sigma_\lambda^2 = \rho$. Hence for any $t > 0$
$$\rho^2 1_{\rho \leq \sqrt{\rho}t} + D a_{j,D}^2 1_{|a_{j,D}| \leq \sqrt{\rho}t} \leq R^{2-4s} t^{4s}.$$

If $|a_{j,D}| \leq \rho$, then it is enough to have
$$\rho^2 + D a_{j,D}^2 \leq R^{2-4s} \rho^{2s} \tag{4.17}$$

and
$$D a_{j,D}^2 \leq R^{2-4s} \left(\frac{a_{j,D}^2}{\rho}\right)^{2s}. \tag{4.18}$$

If the parameters satisfy these equations, then
$$\mathcal{R}(W_s(R) \cap B_{2,\Gamma}^s(R') \cap \mathcal{L}_{1,2,\infty}(R'')) \geq \mathcal{R}(\mathcal{M}_{j,D}),$$

where $\mathcal{R}(W_s(R) \cap B_{2,\Gamma}^s(R') \cap \mathcal{L}_{1,2,\infty}(R''))$ and $\mathcal{R}(\mathcal{M}_{j,D})$ are respectively the minimax risks associated with $W_s(R) \cap B_{2,\Gamma}^s(R') \cap \mathcal{L}_{1,2,\infty}(R'')$ and $\mathcal{M}_{j,D}$. By similar arguments to those of the proof of Theorem 5, one obtains
$$\mathcal{R}(\mathcal{M}_{j,D}) \geq \frac{1}{4} \theta D a_{j,D}^2 \inf_{\hat{f} \in \mathcal{M}_{j,D}} (1 - \inf_{f \in \mathcal{M}_{j,D}} \mathbb{P}(\hat{f} = f)).$$

We now use Lemma 6. Recall that (see (4.13))

$$K(\mathbb{P}_{f'_m}, \mathbb{P}_{f_m}) \leq \frac{2}{\rho} n D a_{j,D}^2.$$

Hence

$$\mathcal{R}(\mathcal{M}_{j,D}) \geq \frac{(1-B)\theta}{4} D a_{j,D}^2$$

as soon as the mean Kullback Leibler distance is small enough, which is implied by

$$\frac{2}{\rho} n D a_{j,D}^2 \leq B\sigma D \log(2^j/D). \tag{4.19}$$

Let us take $j$ such that $2^j \leq n/\log n \leq 2^{j+1}$ and with $D \leq 2^j$,

$$a_{j,D}^2 = \frac{\rho^2}{4n} \log(2^j/D).$$

First note that (4.19) is automatically fulfilled as soon as $\rho \leq 2B\sigma$, that is true if $\rho$ an absolute constant small enough. Then

$$\rho + 2^{j/2}|a_{j,D}| \leq \rho + 2^{j/2}\sqrt{\frac{\rho^2 \log n}{4n}} \leq 1.5\rho.$$

So, if $\rho$ is an absolute constant small enough, (4.15) is satisfied. Moreover

$$2^{1+j/2}|a_{j,D}| \leq 2^{1+j/2}\sqrt{\frac{\rho^2 \log n}{4n}} \leq \rho.$$

This gives (4.14). Now, take an integer $D = D_n$ such that

$$D_n \sim_{n \to \infty} R^{2-4s} \left(\frac{n}{\log n}\right)^{1-2s}.$$

For $n$ large enough, $D_n \leq 2^j$ and $D_n$ is feasible. We have for $R$ fixed,

$$a_{j,D_n}^2 \sim_{n \to \infty} C_s \rho^2 \frac{\log n}{n},$$

where $C_s$ is a constant only depending on $s$. Therefore,

$$\rho + 2^{js}\sqrt{D_n}|a_{j,D_n}| = \rho + \sqrt{C_s}\rho R^{1-2s} + o_n(1).$$

Since $R^{1-2s} \leq R'$ it is sufficient to take $\rho$ small enough but constant depending only on $s$ to obtain (4.16). Moreover,

$$D_n a_{j,D_n}^2 \sim_{n \to \infty} C_s \rho^2 R^{2-4s} \left(\frac{\log n}{n}\right)^{2s}.$$

Hence (4.17) is equivalent to $\rho^2 < R^{2-4s}\rho^{2s}$. Since $R \geq 1$, this is true as soon as $\rho < 1$. Finally (4.18) is equivalent, when $n$ tends to $+\infty$, to

$$C_s\rho^2 \leq (C_s\rho)^{2s}.$$

Once again this is true for $\rho$ small enough depending on $s$. As we can choose $\rho$ not depending on $R, R', R''$, this concludes the proof.

Corollary 1 is completely straightforward once we notice that if $R' \geq R$ then for every $s$, $R' \geq R^{2-4s}$.

**Acknowledgment.** The authors acknowledge the support of the French Agence Nationale de la Recherche (ANR), under grant ATLAS (JCJC06 137446) "From Applications to Theory in Learning and Adaptive Statistics". We would like to warmly thank Lucien Birgé for his advises and his encouragements.